\documentclass[11pt,a4paper]{article}
\usepackage{amsmath}
\usepackage{amsthm}
\usepackage{amssymb}
\usepackage{dsfont}
\usepackage[T1]{fontenc}
\usepackage[french,english,russian]{babel}
\usepackage[utf8]{inputenc}
\usepackage{mathrsfs}
\usepackage{lmodern}
\usepackage{graphicx}
\newtheorem{lemma}{Lemma}
\newtheorem{sublemma}{Sublemma}
\newtheorem{prop}{Proposition}
\newtheorem{theorem}{Theorem}
\newtheorem*{theorems}{Theorem}

\newtheorem{defi}{Definition}
\makeatletter
\def\moverlay{\mathpalette\mov@rlay}
\def\mov@rlay#1#2{\leavevmode\vtop{%
   \baselineskip\z@skip \lineskiplimit-\maxdimen
   \ialign{\hfil$\m@th#1##$\hfil\cr#2\crcr}}}
\newcommand{\charfusion}[3][\mathord]{
    #1{\ifx#1\mathop\vphantom{#2}\fi
        \mathpalette\mov@rlay{#2\cr#3}
      }
    \ifx#1\mathop\expandafter\displaylimits\fi}
\makeatother

\title{Limit laws in the lattice problem. \\
II. The case of ovals.}
\author{Julien Trevisan}
\begin{document}
\selectlanguage{english}
\maketitle
\bigskip
\selectlanguage{french}
\begin{abstract}
Nous étudions l'erreur du nombre de points d'un réseau unimodulaire qui tombent dans une ellipse dilaté et centré sur $0$. On montre d'abord que l'étude de l'erreur, lorsqu'elle est normalisée par $\sqrt{t}$ où $t$ est le paramètre de dilatation de l'ellipse, lorsque ce paramètre tend vers l'infini et lorsque le réseau considéré est aléatoire, se ramène à l'étude d'une transformée de Siegel $\mathcal{S}(f_{t})(L)$ qui dépend de $t$. Ensuite, on se ramène à l'étude du comportement asymptotic d'une transformée de Siegel avec poids aléatoires, $\mathcal{S}(F)(\theta,L)$ où $\theta$ est un second paramètre aléatoire. Enfin, on montre que cette dernière quantité converge presque sûrement et on étudie l'existence des moments de sa loi.
\end{abstract}
\selectlanguage{english}
\begin{abstract}
We study the error of the number of unimodular lattice points that fall into a dilated and centred ellipse around $0$. We first show that the study of the error, when the error is normalized by $\sqrt{t}$ with $t$ the parameter of dilatation of the ellipse, when $t$ tends to infinity and when the lattice is random, is reduced to the study of a Siegel transform $\mathcal{S}(f_{t})(L)$ that depends on $t$. Then, by making $t \rightarrow \infty$, we see that $\mathcal{S}(f_{t})$ converges in law towards a modified Siegel transform with random weights $\mathcal{S}(F)(\theta,L)$ where $\theta$ is a second random parameter. Finally, we show that this last quantity converges almost surely and we study the existence of the moments of its law.
\end{abstract}
\section{Introduction}
The lattice point problem is an open problem in the Geometry of Numbers, at least since Carl Friedrich Gauss took interest in which became the famous Gauss circle problem. The general problem states as followed. \\
Let $d$ be an integer greater than $1$. We recall the following definition : 
\begin{defi}
A subset $L$ of $\mathbb{R}^{d}$ is a lattice if it is a subgroup of $\mathbb{R}^{d}$ such that $L$ is discrete and $\text{span}(L) = \mathbb{R}^{d}$. 
\end{defi}
 Let $P$ be a measurable subset of $\mathbb{R}^{d}$ of non-zero finite Lebesgue measure. We want to evaluate the following cardinal number when $t \rightarrow \infty$ : $$ N(tP + X, L) = | (tP + X) \cap L|$$ where $X \in \mathbb{R}^{d}$, $L$ is a lattice of $\mathbb{R}^{d}$ and $t P + X$ denotes the set $P$ dilated by a factor $t$ relatively to $0$ and then translated by the vector $X$.  \\
Under mild regularity conditions on the set $P$, one can show that : 
$$ N(tP + X, L) = t^{d}\frac{\text{Vol}(P)}{\text{Covol}(L)} + o(t^{d}) $$
where $o(f(t))$ denotes a quantity such that, when divided by $f(t)$, it goes to $0$ when $t \rightarrow \infty$. \\
We are interested in the error term $$\mathcal{R}(tP + X,L) = N(tP + X, L) - t^{d}\frac{\text{Vol}(P)}{\text{Covol}(L)} \textit{.}$$
In the case where $d=2$ and where $P$ is the unit disk $\mathbb{D}^{2}$, Hardy's conjecture in $\cite{hardy1917average}$ stipulates that we should have for all $\epsilon > 0$, $$\mathcal{R}(t \mathbb{D}^{2}, \mathbb{Z}^{2}) = O(t^{\frac{1}{2}+\epsilon}) \textit{.}$$
One of the result in this direction has been established by Iwaniec and Mozzochi in $\cite{iwaniec1988divisor}$. They have proven that for all $\epsilon > 0$, $$\mathcal{R}(t \mathbb{D}^{2}, \mathbb{Z}^{2}) = O(t^{\frac{7}{11}+\epsilon}) \textit{.}$$
This result has been recently improved by Huxley in \cite{huxley2003exponential}. Indeed, he has proven that : $$\mathcal{R}(t \mathbb{D}^{2}, \mathbb{Z}^{2}) = O(t^{K} \log(t)^{\Lambda}) $$
where $K = \frac{131}{208} $ and $\Lambda = \frac{18627}{8320}$.  \\
In dimension 3, Heath-Brown has proven in $\cite{heath2012lattice}$ that : 
$$\mathcal{R}(t \mathbb{D}^{3}, \mathbb{Z}^{3}) = O(t^{\frac{21}{16}+\epsilon}) \textit{.}$$
These last two results are all based on estimating what are called $\textit{exponential sums}$.\\ Furthermore, in both cases, the error is considered in a deterministic way. Another approach was followed by Bleher, Cheng, Dyson and Lebowitz in $\cite{bleher1993distribution}$. They took interest in the case where the dilatation parameter $t$ is random. More precisely, they took interest in the case where $t$ is being distributed according to the measure $\rho(\frac{t}{T}) dt$ (that is absolutely continuous relatively to Lebesgue measure) and where $\rho$ is a probability density on $[0,1]$ and $T$ is parameter that goes to infinity. In that case, they showed the following result : 
\begin{theorems}[$\cite{bleher1993distribution}$]
There exists a probability density $p$ on $\mathbb{R}$ such that for every piecewise continuous and bounded function $g : \mathbb{R} \longrightarrow \mathbb{R}$, 
$$ \lim_{T \rightarrow \infty} \frac{1}{T} \int_{0}^{T} g(\frac{\mathcal{R}(t \mathbb{D}^{2}, \mathbb{Z}^{2})}{\sqrt{t}}) \rho(\frac{t}{T}) dt = \int_{\mathbb{R}} g(x) p(x) dx \textit{.}$$
Furthermore $p$ can be extended as an analytic function over $\mathbb{C}$ and verifies that for every $\epsilon > 0$, $$p(x) = O(e^{-|x|^{4- \epsilon}})$$ when $x \in \mathbb{R}$ and when $|x| \rightarrow \infty$. 
\end{theorems}
Since the works of Kesten $\cite{kesten60u}$ and $\cite{kesten62u}$, a lot of work has been done in estimating errors terms : we can cite, for example, $\cite{gorodnikcounting}$ and $\cite{dolgopyat2014deviations}$. Following their approaches, we have taken interest in what is happening to $\mathcal{R}$ when $L$ is a random unimodular lattice, when $P=\mathcal{E}$ is an ellipse centred on $0$ and when $t$ is a deterministic parameter that goes to infinity. \\
Let's recall the following definition : 
\begin{defi}
The coarea of a lattice $L$ of $\mathbb{R}^{2}$ is the Lebesgue measure of a measurable fundamental set. Furthermore, a lattice is said to be unimodular if its corea is equal to $1$.
\end{defi}
We denote by $\mathscr{S}_{2}$ the space of unimodular lattices and it can be seen as the quotient space $SL_{2}(\mathbb{R})/SL_{2}(\mathbb{Z})$. We denote by $\mu_{2}$ the unique Haar probability measure on it. Let us set : 
$$\Pi = \{ (k_{1},k_{2}) \in \mathbb{Z}^{2} \text{ } | \text{ } k_{1} \wedge k_{2} = 1 \textit{, } k_{1} \geqslant 0  \} $$
where we agree that if $k_{1} = 0$, $k_{1} \wedge k_{2} = 1$ means that $k_{2} = 1$. \\
Let's set : $$\mathbb{T}^{\infty} = (\mathbb{T}^{1})^{\Pi}$$ 
where $\mathbb{T}^{1} = \mathbb{R}/\mathbb{Z}$ and let's call $\lambda_{\infty}$ the normalized Lebesgue measure product over $\mathbb{T}^{\infty}$. \\
Let's recall other usual definitions : 
\begin{defi}
\label{def3}
For a lattice $L \in \mathscr{S}_{2}$, we say that a vector of $L$ is prime if it not a non-trivial integer multiple of another vector of $L$. In fact, a vector $ l \in L$ is prime if, and only if, $M l \in M L$ with $M \in SL_{2}(\mathbb{R})$ is prime and a vector $(k_{1},k_{2}) \in \mathbb{Z}^{2}$ is prime if, and only if, $k_{1} \wedge k_{2} = 1$. \\
\end{defi}

\begin{defi}
\label{def1}
For every $ i \in \{ 1, 2 \}$, we call $$\lVert L \rVert _{i} = \min \{ r > 0 \text{ } | \text{ } B_{f}(0,r) \text{ contains } i \text{ vectors of L lineary independent} \}  $$
where $B_{f}(0,r)$ is the closed centred ball on $0$ for the norm $\lVert \cdot \rVert$ of radius $r$. In fact, for almost all $L \in \mathscr{S}_{2}$, $\lVert L \rVert_{2} > \lVert L \rVert_{1} $ and there exists only one couple of vectors $(e_{1}(L),e_{2}(L))$ such that $(e_{1}(L))_{1} > 0$, $ \lVert e_{1}(L) \rVert = \lVert L \rVert_{1}$, $(e_{2}(L))_{1} > 0$ and $\lVert e_{2}(L) \rVert = \lVert L \rVert_{2}$. \\
\end{defi}
With these notations, one has that $ e \in L$ is prime if, and only if, $e$ can be written as $ e=k_{1} e_{1}(L) + k_{2} e_{2}(L)$ with $k_{1} \wedge k_{2} = 1$. \\
\begin{defi}
We call $\mathcal{P}_{+}(L)$ the set of vectors $e$ of $L$ such that $k_{1} \geqslant 0$ where $e=k_{1} e_{1}(L) + k_{2} e_{2}(L)$ and where $k_{1} \wedge k_{2} = 1$. It is a set that is well-defined for almost all, relatively to $\mu_{2}$, $L \in \mathscr{S}_{2}$ and all the vectors of $\mathcal{P}_{+}(L)$ are prime vectors according to the previous remark. \\
\end{defi}
\begin{defi}
We say that a real random variable $Z$ is symmetrical if $$\mathbb{P}_{Z} = \mathbb{P}_{-Z} \textit{.}$$
\end{defi}
Let $\tilde{\mu}_{2}$ be a probability measure that has a smooth bounded density $\sigma$ with respect to $\mu_{2}$. \\
Let $\mathcal{E}$ be an ellipse centred around $0$. \\
We denote by $\lVert \cdot \rVert$ the usual euclidean norm over $\mathbb{R}^{2}$.\\
The main result of this article is the following theorem : 
\begin{theorem}
\label{thm2}
For every real numbers $a <b$, 
$$\lim_{t \rightarrow \infty} \tilde{\mu}_{2} \left( \frac{\mathcal{R}(t \mathcal{E}, L)}{\sqrt{t}} \in [a,b] \right) = (\lambda_{\infty} \times \tilde{\mu}_{2} ) \left(  S(\theta,L) \in [a,b] \right) $$
where $\theta = (\theta_{e}) \in \mathbb{T}^{\infty}$,
$$S(\theta,L) = \frac{2}{\pi} \sum_{e \in  \mathcal{P}_{+}(L) } \frac{\phi(\theta_{e})}{\lVert e \rVert^{\frac{3}{2}}}  $$ 
with 
\begin{equation}
\label{eq501}
\phi(\theta_{e})= \sum_{m \geqslant 1} \frac{\cos(2 \pi m \theta_{e} - \frac{3 \pi}{4})}{m^{\frac{3}{2}}} \textit{.}
\end{equation}
Furthermore, $S(\theta,L)$ : \begin{itemize}
\item converges almost surely
\item is symmetrical and its expectation is equal to $0$
\item admits a moment of order $1+\kappa$ for any $0 \leqslant \kappa < \frac{1}{3}$ 
\item does not admit a moment of order $\frac{4}{3}$ when $\sigma(L) \geqslant m$ where $m > 0$ and where $L$ belongs to an event of the form $( \lVert L \rVert < \alpha )$ with $\alpha > 0$. 
\item when there exists $\alpha > 0 $ such that $\tilde{\mu}_{2}(\{ L \in \mathscr{S}_{2} \text{ } | \text{ } \lVert L \rVert_{1} < \alpha \}) = 0$ then $S(\theta,L)$ admits moments of all order $1 \leqslant p < \infty$. 
\end{itemize}
\end{theorem}
The proof of Theorem $\ref{thm2}$ can be reduced to the case $\mathcal{E} = \mathbb{D}^{2}$. Indeed, since the parameter of dilatation goes to infinity, we can replace it by $\frac{t}{\text{Area}(\mathcal{E})}$ and so the study is reduced to the case where $\mathcal{E}$ is of area $1$. Furthermore, we note that in that case there exists $M \in SL_{2}(\mathbb{R})$ such that $M \mathcal{E} = \mathbb{D}^{2}$ and so the study is reduced to the case where $P = \mathbb{D}^{2}$ by considering that $L$ is distributed according to the probability measure $(M^{-1})_{*} \tilde{\mu}_{2}$. As a consequence, for the rest of the article, we are going to suppose $$\mathcal{E} = \mathbb{D}^{2} \textit{.}$$
We notice that, in the case where $\sigma(L) \geqslant m$ where $m > 0$ and where $L$ belongs to an event of the form $( \lVert L \rVert < \alpha )$ with $\alpha > 0$, the tail of the distribution is by far bigger than the tail of the limit distribution of $\cite{bleher1993distribution}$. Indeed, in our case, the limit distribution admits of order less than $\frac{4}{3}$ but does not admit a moment of order larger than $\frac{4}{3}$. In the case of $\cite{bleher1993distribution}$, the limit distribution admits moment of all order $1 \leqslant p < \infty$. In fact, they have better : the decreasing of the tail is at least estimated by $\int_{x \geqslant \alpha} e^{-|x|^{4-\epsilon}} $ where $\alpha \rightarrow \infty$ and $\epsilon$ is an arbitrary positive parameter.
\\
\\
The next section is a preliminary dedicated to a summation formula. It will enable us to give a heuristic explanation of the approach that we will follow. At the end of the section, we will give the plan of the paper.
\section{Heuristic explanation and plan of the proof}
We are first going to recall a summation formula (about this subject, see also section 2 of $\cite{bleher1992distribution}$). This summation formula is based on the Poisson formula. Yet it does not apply in the case treated by Bleher, nor in our case. In fact, it applies after having $\textit{regularized}$ the problem (see the beginning of section 3 for more details). So it will give us a good insight of the problem and that is the reason why we are going to recall this summation formula. \\
After doing so, we will apply it to our case. It will finally enable us to give heuristic explanations and a detailed plan of the proof of Theorem $\ref{thm2}$. \\
Let $\gamma$ be a simple, smooth, closed, convex and with a positive curvature curve in $\mathbb{R}^{2}$ such that $(0,0)$ lies inside $\gamma$. For $\Omega_{\gamma}$ the domain enclosed by $\gamma$, we call
$ \mathbf{1}_{t \Omega_{\gamma}} $
the characteristic function of $t \Omega_{\gamma}$. \\
Then one has that : $$N(t \Omega_{\gamma}, \mathbb{Z}^{2}) = \sum_{n \in \mathbb{Z}^{2}} \mathbf{1}_{t \Omega_{\gamma}} (n) \textit{.} $$ 
By the Poisson summation formula (yet, as we said earlier, it does not apply here because $\mathbf{1}_{t \Omega_{\gamma}}$ is not smooth enough and so there is a convergence problem), one has that : 
$$N(t \Omega_{\gamma}, \mathbb{Z}^{2}) = \sum_{n \in \mathbb{Z}^{2}} \widetilde{\mathbf{1}_{t \Omega_{\gamma}}} (2 \pi n) $$ 
with $  \widetilde{\mathbf{1}_{t \Omega_{\gamma}}}$ the Fourier transform of $\mathbf{1}_{t \Omega_{\gamma}}$ defined by $$\widetilde{\mathbf{1}_{t \Omega_{\gamma}}} (\xi) = \int_{\mathbb{R}^{2}} e^{i x \xi}  \mathbf{1}_{t \Omega_{\gamma}}(x) dx \textit{.}$$
Note that $\widetilde{\mathbf{1}_{t \Omega_{\gamma}}} (0) = \text{Area}(t \Omega_{\gamma})$. So we obtain that : 
\begin{equation}
\label{eq1}
\mathcal{R}(t \Omega_{\gamma},\mathbb{Z}^{2}) = \sum_{n \in \mathbb{Z}^{2}-\{ 0 \}}  \widetilde{\mathbf{1}_{t \Omega_{\gamma}}} (2 \pi n) \textit{.} 
\end{equation}
We need some notations before moving forward. For $\xi \in \mathbb{R}^{2}- \{ 0 \}$, let's call $x_{\gamma}(\xi)$ the point on $\gamma$ where the outer normal to $\gamma$ coincides with $\frac{\xi}{\lVert \xi \rVert}$ where $\lVert \xi \rVert$ refers to the usual euclidean norm of $\mathbb{R}^{2}$. Let's call also $\rho_{\gamma}(\xi)$ denote the curvature radius of $\gamma$ at $x_{\gamma}(\xi)$ and $$ Y_{\gamma}(\xi) = < \xi, x_{\gamma}(\xi)> $$
where $<\cdot,\cdot>$ denotes the usual scalar product over $\mathbb{R}^{2}$. \\
With these notations, Bleher in $\cite{bleher1992distribution}$ obtains heuristically by giving an asymptotic expression of $ \widetilde{\mathbf{1}_{t \Omega_{\gamma}}}$ and by using ($\ref{eq1}$)  : 
\begin{equation}
\label{eq2}
\frac{\mathcal{R}(t \Omega_{\gamma},\mathbb{Z}^{2})}{\sqrt{t}} = \frac{1}{\pi} \sum_{n \in \mathbb{Z}^{2}-\{ 0 \}} \frac{\sqrt{\rho_{\gamma}(n)} \cos(2 \pi t Y_{\gamma} (n) - \frac{3 \pi}{4})}{\lVert n \rVert^{\frac{3}{2}}} + O(\frac{1}{t})
\end{equation}
where $Y=O(f(t))$ is a quantity such that $$|Y| \leqslant C |f(t)| $$ with $C \geqslant 0$. \\
Before saying what this summation formula gives in our case we just need to recall the definition of a dual lattice : 
\begin{defi}
\label{def2}
For $L \in \mathscr{S}_{2}$, its dual lattice $L^{\perp}$ is defined in the following way
$$L^{\perp} = \{ x \in \mathbb{R}^{2} \text{ } | \text{ } \forall l \in L \text{, } <l,x> \in \mathbb{Z} \} \textit{.} $$
Furthermore, if $M$ is a matrix in $SL_{2}(\mathbb{R})$ that represents $L$ then $(M^{-1})^{T}$ is a matrix that represents $L^{\perp}$ where $(M^{-1})^{T}$ is the transpose matrix of $M^{-1}$.
\end{defi}
With the previously exposed summation formula, in our case we get, after some calculations (that we will expose later) : 
$$\frac{\mathcal{R}(t \mathbb{D}^{2},L)}{\sqrt{t}} = \frac{1}{\pi} \sum_{l \in L^{\perp}-\{ 0 \} } \frac{ \cos(2 \pi t \lVert l \rVert - \frac{3 \pi}{4})}{\lVert l \rVert^{\frac{3}{2}}} + O(\frac{1}{t}) $$ 
where $L^{\perp}$ is the dual lattice of $L$.
Let's recall the following definition : 
\begin{defi}
\label{def10}
Let $f$ be a function on $\mathbb{R}^{2}$.  \\
Let's define formally : 
$$ \mathcal{S}(f)(L) = \sum_{\substack{e \in L \\ e \textit{ prime} }} f(e)  \textit{.}$$
$\mathcal{S}(f)$ is called the $\textit{Siegel}$ transform of $f$. \\
If $f$ has a compact support, $\mathcal{S}(f)$ is well-defined.
\end{defi}
Let's set $f_{t}(x) = \frac{1}{\pi} \frac{ \cos(2 \pi t \lVert x \rVert - \frac{3 \pi}{4})}{\lVert x \rVert^{\frac{3}{2}}} \mathbf{1}_{\mathbb{R}^{2}-\{0 \}}(x) \textit{.}$
As a consequence of the summation formula, heuristically the asymptotic behaviour of $\frac{\mathcal{R}(t \mathbb{D}^{2},L)}{\sqrt{t}}$ is the same as the asymptotic behaviour of the following Siegel transform
\begin{equation} 
\label{eq500}
\mathcal{S}(f_{t})(L^{\perp}) = \frac{1}{\pi} \sum_{l \in L^{\perp}-\{ 0 \} } \frac{ \cos(2 \pi t \lVert l \rVert - \frac{3 \pi}{4})}{\lVert l \rVert^{\frac{3}{2}}} \textit{.} 
\end{equation}
We can, and we will, study it when $t \rightarrow \infty$ by replacing $L^{\perp}$ by $L$ (because $L$ is supposed to be random). \\
By the way, here we see that the randomness on $t$ in the case treated in $\cite{bleher1993distribution}$ is somewhat replaced by the randomness on $ \lVert l \rVert$ where $l \in L-\{ 0 \}$ with $L$ random and distributed according to $\tilde{\mu}_{2}$. Yet, in $\cite{bleher1993distribution}$ and in $\cite{bleher1992distribution}$, this randomness on $t$ enables the authors to use results about what is called $\textit{almost periodic}$ function (developed in $\cite{besicovitch1954almost}$ and in $  \cite{heath1992distribution} $). More precisely, the authors use the following theorem : for every $F : \mathbb{R}_{+} \longrightarrow \mathbb{R}$, if for every $\epsilon > 0$, there exists a trigonometric polynomial $P_{\epsilon}$ such that 
$$ \limsup_{T \rightarrow \infty} \frac{1}{T} \int_{0}^{T} \min(1,|F(t)-P_{\epsilon}(t)|) \leqslant \epsilon $$ then $F(t)$ admits an asymptotic distribution.\\ Here, in our case, we can not work in this framework (because of the difference of nature of the randomness). So, we have to take another approach and the approach we are going to follow is directly taken off $\cite{dolgopyat2014deviations}$. \\
We now give the main steps of the proof of Theorem $\ref{thm2}$.\\
\\
$\textit{First step.}$ By regularizing the problem, we are going to show that the quantity $S_{A,prime}(L,t)$, for $A > 0$ a fixed parameter that is taken large enough, when $t \rightarrow \infty$, is close, in probability, to $\frac{\mathcal{R}(t \mathbb{D}^{2},L)}{\sqrt{t}}$ where $S_{A,prime}(L,t)$ is defined by : 
\begin{equation}
\label{eq44}
S_{A,prime}(L,t) = \frac{1}{\pi} \sum_{\substack{ l \in L^{\perp} \textit{ prime}  \\ 0 < \lVert l \rVert \leqslant A }} \frac{ 1}{\lVert l \rVert^{\frac{3}{2}}} \sum_{m \in \mathbb{N}-\{ 0 \}} \frac{ \cos(2 \pi t m \lVert l \rVert - \frac{3 \pi}{4})}{m^{\frac{3}{2}}} \textit{.}
\end{equation}
It is different from the right-hand side of equation ($\ref{eq500}$) because, first, the sum has been cut and it has to be cut because of a problem of convergence of the right-hand side of equation ($\ref{eq500}$) ; second, we take into account a phenomenon of multiplicity (if $l$ appears in the sum, $2l$ is also going to appear). Bleher in $\cite{bleher1993distribution}$ did also take into account such a phenomenon. He did it in order to get independence at infinity, as it was also done in $\cite{bassam}$, and we do that for the same goal.  \\
By using a simple parity argument, by using the remark that follows Definition $\ref{def3}$ and by replacing $L$ by $L^{\perp}$, one has that : 
\begin{equation}
\label{eq200}
S_{A,prime}(L^{\perp},t) = \frac{2}{\pi} \sum_{\substack {k_{1} \wedge k_{2} = 1 \\ k_{1} \geqslant 0 \\ \lVert k_{1} e_{1}(L) + k_{2} e_{2}(L) \rVert \leqslant A }} \frac{\phi(t \lVert k_{1} e_{1}(L) + k_{2} e_{2}(L) \rVert)}{\lVert k_{1} e_{1}(L) + k_{2} e_{2}(L) \rVert^{\frac{3}{2}}} 
\end{equation}
where the function $\phi$ was defined by the equation $(\ref{eq501})$. We have done that so from this stage onwards we consider vectors of $\mathcal{P}_{+}(L)$ with a fixed indexation (that does not depend on $L$). Furthermore, this indexation will be very useful for the second step (for more details, see section 4). \\
\\
$\textit{Second step.}$ 
We will show that the family of variables $(t  \lVert k_{1} e_{1}(L) + k_{2} e_{2}(L) \rVert)$, whose values are in $\mathbb{R}/ \mathbb{Z}$, are, when $t \rightarrow \infty$, independent from one another and are identically distributed according to the normalized Haar measure on $\mathbb{R}/ \mathbb{Z}$. The idea here is basically the same as in $\cite{bleher1993distribution}$ and in $\cite{heath1992distribution}$ where the respective authors used the fact that the square roots of square free integers are $\mathbb{Z}$-free. In our case, to prove the result, we will decompose the space of unimodular lattices into small geodesic segments, calculate the Taylor series of $\lVert k_{1} e_{1}(L) + k_{2} e_{2}(L) \rVert$ at order 1 on such a segment and show that the coefficients of order 1 are $\mathbb{Z}$-free. \\
We will also prove that these variables are independent, when $t \rightarrow \infty$, from the variables $L$ due to the presence of the factor $t$. \\
\\
$\textit{Third step.}$
Thanks to the first and second step, we will see that the asymptotic distribution of $\frac{\mathcal{R}(t \mathbb{D}^{2},L)}{\sqrt{t}}$ is the distribution of
\begin{equation}
\label{eq302}
 S(\theta,L) = \sum_{\substack {k = (k_{1},k_{2}) \\ k_{1} \wedge k_{2} = 1 \\ k_{1} \geqslant 0 }} \frac{\phi(\theta_{k})}{\lVert k_{1} e_{1}(L) + k_{2} e_{2}(L) \rVert^{\frac{3}{2}}}  
 \end{equation}
under the assumption that the quantity $S_(\theta,L)$ is well-defined and where $\theta = (\theta_{k}) \in \mathbb{T}^{\infty}$ is distributed according to $\lambda_{\infty}$, $ L \in \mathscr{S}_{2}$ is distributed according to $\tilde{\mu}_{2}$  It should be noted that the $\left( \phi(\theta_{k}) \right)$ is a family of independent and identically distributed real random variables. Furthermore, all of them are symmetrical and have a compact support. The third step of the proof will consist in showing a general result. Let the $Z_{k}$ be real non-zero independent identically distributed real random variables from a probability space $\Omega \ni \omega$ that are symmetrical and have a compact support.  Then we are going to show that
\begin{equation}
\label{eq301}
\tilde{S}_{A}(\omega,L) = \sum_{\substack {k = (k_{1},k_{2}) \\ k_{1} \wedge k_{2} = 1 \\ k_{1} \geqslant 0 \\ \lVert k_{1} e_{1}(L) + k_{2} e_{2}(L) \rVert \leqslant A}} \frac{Z_{k}(\omega)}{\lVert k_{1} e_{1}(L) + k_{2} e_{2}(L) \rVert^{\frac{3}{2}}}
\end{equation}
converges almost surely when $A \rightarrow \infty$. $\tilde{S}$ looks like a Siegel transform but it is not because of the numerator of the considered terms. We are going to call it a modified Siegel transform with random weights. An object of this type was already studied in $\cite{dolgopyat2014deviations}$. As a consequence, $S_(\theta,L)$, as the limit when $A \rightarrow \infty$ of a particular $\tilde{S}$, will be well-defined. \\
We will also study the existence of moments of the almost sure limit. In particular, we are going to see that the optimal $\kappa$ is $\frac{1}{3}$ (see statement of Theorem $\ref{thm2}$). Heuristically, it is because the magnitude of $\tilde{S}$ is given by the term $$ \frac{Z_{1,0}}{\lVert L \rVert_{1}^{\frac{3}{2}}}$$ with $Z_{1,0} \neq 0$.\\
After doing all of that, we will finally get the validity of Theorem $\ref{thm2}$. \\
\\
$\textbf{Plan of the paper.}$
The next section will be dedicated to deal with the first step of the proof, namely it will show that $\frac{\mathcal{R}(t \mathbb{D}^{2},L)}{\sqrt{t}}$ is close in probability with $S_{A,prime}(L^{\perp},t) $
when $A$ is a fixed parameter taken large enough and $t$ goes to infinity (see Proposition $\ref{prop8}$). We have to "cut" the sum because of the problem of convergence of the Fourier series of $X  \longmapsto \frac{\mathcal{R}(t \mathbb{D}^{2}+X,L)}{\sqrt{t}}$ which is due to the lack of regularity of the indicator function $\mathbf{1}_{t \mathbb{D}^{2}}$. To prove this, we are going to proceed by $\textit{regularization}$ which means here that we are going to smooth the indicator function $\mathbf{1}_{t \mathbb{D}^{2}}$ via a Gaussian kernel. \\
In section 4, we tackle the second step of the proof, that is the fact that we have independence when $t \rightarrow \infty$ of the $(t \lVert k_{1}e_{1}(L) + k_{2} e_{2}(L) \rVert)$, that they are identically distributed according to the normalized Haar measure over $\mathbb{R}/\mathbb{Z}$ independent from $L$ and so that, when $t \rightarrow \infty$, $\frac{\mathcal{R}(t \mathbb{D}^{2},L)}{\sqrt{t}}$ has the same distribution of $S(\theta,L)$. \\
In section 5, which is the last section, we are going to tackle with the third step of the proof, namely study the convergence of $\tilde{S}(\omega, L, A)$ when $A \rightarrow \infty$ and the existence of moments of its limit. \\
In the rest of the article, all the calculus of expectation $\mathbb{E}$, of variance $\text{Var}$ and of probability $\mathbb{P}$ will be made according to the measure $\tilde{\mu}_{2}$. Furthermore, the expression typical is going to signify $\tilde{\mu}_{2}-\textit{almost surely}$.
\section{Reduction to the study of the Siegel transform}
The main object of this section is to show the following proposition : 
\begin{prop} 
\label{prop8}
For every $\alpha > 0$, for every $A > 0$ large enough, for every $t$ large enough, one has that :
$$\mathbb{P}(\Delta_{A,prime}(L,t)  \geqslant \alpha) \leqslant \alpha$$
where 
\begin{equation}
\label{eq201}
\Delta_{A,prime}(L,t) = |\frac{\mathcal{R}(t \mathbb{D}^{2},L)}{\sqrt{t}} - S_{A,prime}(L,t)| \textit{ .}
\end{equation}

\end{prop}
This proposition basically says that we can reduce the asymptotical study of $\frac{\mathcal{R}(t \mathbb{D}^{2},L)}{\sqrt{t}}$ to the study of its Fourier transform, taking into account a phenomenon of multiplicity. \\
In fact, due to the triangle inequality, we only have to prove the following two lemmas : 
\begin{lemma}
\label{lemme30}
For every $\alpha > 0$, for every $A > 0$ large enough, for every $t$ large enough, one has that :
$$\mathbb{P}(\Delta_{A}(L,t) \geqslant \alpha) \leqslant \alpha$$
where 
\begin{equation}
\label{eq202}
\Delta_{A}(L,t) = |\frac{\mathcal{R}(t \mathbb{D}^{2},L)}{\sqrt{t}} - H_{A}(L,t)| \textit{ with}
\end{equation}
\begin{equation}
\label{eq10}
H_{A}(L,t) = \frac{1}{\pi} \sum_{\substack{ l \in L^{\perp} \\ 0 < \lVert l \rVert \leqslant A }} \frac{ \cos(2 \pi t \lVert l \rVert - \frac{3 \pi}{4})}{\lVert l \rVert^{\frac{3}{2}}} \textit{.} 
\end{equation}
\end{lemma}
\begin{lemma}
\label{lemme31}
For every $\alpha > 0$, for every $A > 0$ large enough, for every $t$ large enough, one has that :
$$\mathbb{P} (|S_{A,prime}(L,t) - H_{A}(L,t)| \geqslant \alpha) \leqslant \alpha \textit{.}$$
\end{lemma}
\begin{proof}[Proof of Proposition $\ref{prop8}$]
One has that : 
\begin{equation}
\label{eq54}
\Delta_{A,prime}(L,t) \leqslant \Delta_{A}(L,t) + |S_{A,prime}(L,t) - H_{A}(L,t)| \textit{.}
\end{equation}
The Lemma $\ref{lemme30}$ and the Lemma $\ref{lemme31}$ imply then the wanted result.
\end{proof}
Let's say a few words about Lemmas $\ref{lemme30}$ and $\ref{lemme31}$ before following with their respective proofs. The Lemma $\ref{lemme30}$ says that the study of $\frac{\mathcal{R}(t \mathbb{D}^{2},L)}{\sqrt{t}}$ can be reduced to the study of its Fourier transform. The Lemma $\ref{lemme31}$ says that the phenomenon of multiplicity (the fact that for a prime vector $l$, $2 l$, $3 l$ etc. appear in the sum $H_{A}(L,t)$ when $A \rightarrow \infty$) is not so important. We only have to gather all the multiples of a prime vector (which corresponds to the infinite sum over $m$, see equation ($\ref{eq44}$)), so that we focus on prime vectors. 
\subsection{Proof of Lemma $\ref{lemme30}$}
First, we are going to prove the Lemma $\ref{lemme30}$. To do so, we are following closely the approach of $\cite{bleher1992distribution}$, yet with some differences because in our case it is not the radius of dilatation that is random but the lattice (or, equivalently and in a certain sense, the oval). \\
For $x \in \mathbb{R}^{2}$ and $t > 0$, let's define $$\lambda(x ; t) = \frac{t^{2}}{4 \pi} e^{-\frac{t^{2}}{4 \pi} \lVert x \rVert^{2}}$$
and, for $M \in SL_{2}(\mathbb{R})$,
\begin{equation}
\label{eq11}
\lambda_{M}(x;t)=\lambda(M x ; t) \textit{.}
\end{equation}
We recall that : 
\begin{equation}
\label{eq12}
\int_{\mathbb{R}^{2}} \lambda_{M}(x;t) dx = 1 
\end{equation}
and that the Fourier transform of $\lambda_{M}(\cdot; t)$ can be expressed as 
\begin{equation}
\label{eq13}
 \widetilde{\lambda_{M}}(\xi ; t) = e^{- \frac{\lVert (M^{-1})^{T} \xi \rVert^{2}}{t^{2}}} \textit{.}
\end{equation}
We introduce the following function : 
\begin{equation}
\label{eq14}
\chi_{\mathbb{S}^{1}, M}(x;t) = (\mathbf{1}_{t \Omega_{M^{-1}\mathbb{S}^{1}}}*\lambda_{M}(\cdot; t))(x) = \int_{\mathbb{R}^{2}} \mathbf{1}_{t \Omega_{M^{-1}\mathbb{S}^{1}}}(y) \lambda_{M}(x-y; t) dy 
\end{equation}
(it is a regularization of the function $\chi_{t \Omega_{M^{-1}\mathbb{S}^{1}}}$). \\
Let us also set : 
\begin{equation}
\label{eq15}
N_{reg}(t \mathbb{D}^{2}, M) = \sum_{n \in \mathbb{Z}^{2}} \chi_{\mathbb{S}^{1}, M}(n;t) \textit{ and }
\end{equation}
(the index "reg" stands for regularized) 
\begin{equation}
\label{eq16}
F(M,t) = \frac{N_{reg}(t \mathbb{D}^{2}, M) - \text{Area}(t \mathbb{D}^{2})}{\sqrt{t}} \textit{.}
\end{equation}
Let $L$ be a unimodular lattice such that $e_{1}(L)$ and $e_{2}(L)$ are well-defined and let 
\begin{equation}
\label{eq17}
M=[e_{1}(L),e_{2}(L)] \textit{ if } \text{det}([e_{1}(L),e_{2}(L)]) > 0 
\end{equation}
and
\begin{equation}
\label{eq18}
M=[e_{2}(L),e_{1}(L)] \textit{ if } \text{det}([e_{2}(L),e_{1}(L)]) > 0 \textit{.}
\end{equation}
Then $M$ is a matrix that represents $L$. \\
Now, let's call : 
\begin{equation}
\label{eq19}
\Delta_{1}(L,t) = | \frac{\mathcal{R}(t \mathbb{D}^{2},L)}{\sqrt{t}} - F(M,t) | \textit{ and } (\Delta_{2})_{A}(L,t) = | F(M,t) - H_{A}(L,t) |
\end{equation}
so one has that : 
\begin{equation}
\label{eq20}
\Delta_{A}(L,t) \leqslant \Delta_{1}(L,t) + (\Delta_{2})_{A}(L,t)\textit{.}
\end{equation}
The proof of Lemma $\ref{lemme30}$ lies on the two following lemmas : 
\begin{lemma}
\label{lemme32}
$\Delta_{1}(L,t)$ converges almost surely to $0$ when $t \rightarrow \infty$. 
\end{lemma}
\begin{lemma}
\label{lemme33}
For all $\alpha > 0$, for all $A$ large enough, for all $t$ large enough,
$$\mathbb{P}((\Delta_{2})_{A}(L,t) \geqslant \alpha) \leqslant \alpha \textit{.}$$
\end{lemma}
\begin{proof}[Proof of Lemma $\ref{lemme30}$]
It is the direct consequence of equation ($\ref{eq20}$) and of Lemmas $\ref{lemme32}$ and $\ref{lemme33}$.
\end{proof}
Lemma $\ref{lemme32}$ basically tells us that the study of $\frac{\mathcal{R}(t \mathbb{D}^{2},L)}{\sqrt{t}}$ can be reduced to the study of one of its regularized Fourier series, whereas Lemma $\ref{lemme33}$ means that the asymptotical study of this regularized Fourier series can be brought back to the study of the non-regularized Fourier series. \\
The next subsubsection is dedicated to the proof of Lemma $\ref{lemme32}$ and the subsubsection after it is dedicated to the proof of Lemma $\ref{lemme33}$. 
\subsubsection{Proof of Lemma $\ref{lemme32}$}
The proof of Lemma $\ref{lemme32}$ is based on two sublemmas. The first one is the following :
\begin{sublemma}
\label{lemme3}
For all $x \in \mathbb{R}^{2}$, for all $t > 0$, $$|\chi_{\mathbb{S}^{1}, M}(x;t) - \mathbf{1}_{t \Omega_{M^{-1}\mathbb{S}^{1}}}(x)| \leqslant e^{- \frac{t^{2}}{4} \text{dist}(M x,t \mathbb{S}^{1})^{2}} $$ 
where for all $z \in \mathbb{R}^{2}$, $$\text{dist}(z,t \mathbb{S}^{1})= \inf_{y \in t \mathbb{S}^{1}} |z-y| \textit{.}$$
\end{sublemma}
\begin{proof}
One has that : 
$$|\chi_{\mathbb{S}^{1}, M}(x;t) - \mathbf{1}_{t \Omega_{M^{-1}\mathbb{S}^{1}}}(x)|=| \int_{y \notin t \mathbb{D}^{2}} \frac{t^{2}}{4 \pi} e^{- \frac{\lVert Mx - y \rVert ^{2}}{4}} dy| \textit{ if } Mx \in t \mathbb{D}^{2}$$
and 
$$|\chi_{\mathbb{S}^{1}, M}(x;t) - \mathbf{1}_{t \Omega_{M^{-1}\mathbb{S}^{1}}}(x)|=| \int_{y \in t \mathbb{D}^{2}} \frac{t^{2}}{4 \pi} e^{- \frac{\lVert Mx - y \rVert ^{2}}{4}} dy| \textit{ if } Mx \notin t \mathbb{D}^{2}$$
because of equation $(\ref{eq12})$ and by making the change of variable $y = Mu$. \\
The proof of Lemma 3.2 from $\cite{bleher1992distribution}$ gives the wanted result.
\end{proof}
The second sublemma gives an estimate of $dist(Mn,t \mathbb{S}^{1})$. It states in the following way : 
\begin{sublemma}
\label{lemme4}
For almost all $L \in \mathscr{S}_{2}$, for all $t$ large enough, for all $n \in \mathbb{Z}^{2}-\{ 0 \}$, we have that : 
$$dist(Mn,t \mathbb{S}^{1}) \leqslant C t \lVert n \rVert $$
where $C = \frac{(\lVert e_{1}(L) \rVert + \lVert e_{2}(L) \rVert)}{\lVert L \rVert_{1}} \textit{.}$
\end{sublemma}
We recall the reader that $\lVert L \rVert_{1}$ was defined in Definition $\ref{def1}$.
\begin{proof}
By definition of $\text{dist}$, we have that : 
\begin{equation}
\label{eq21}
dist(Mn,t \mathbb{S}^{1}) \leqslant \lVert Mn - \frac{M n}{\lVert Mn \rVert}t \rVert = | 1 - \frac{t}{\lVert Mn \rVert} | \lVert Mn \rVert  \textit{.}
\end{equation}
Furthermore, one has that : 
\begin{equation}
\label{eq22}
\lVert Mn \rVert \leqslant \lVert M \rVert \lVert n \rVert
\end{equation}
where $\lVert M \rVert$ is the matrix norm subordinated to $\lVert \cdot \rVert$. \\
With the equation $(\ref{eq17})$ and equation $(\ref{eq18})$, one can see that : 
\begin{equation}
\label{eq23} 
\lVert M \rVert \leqslant \lVert e_{1}(L) \rVert + \lVert e_{2}(L) \rVert \textit{.}
\end{equation}
Also, we have that : 
\begin{equation}
\label{eq24}
| 1 - \frac{t}{\lVert M n \rVert} | \leqslant \frac{t}{\lVert L \rVert_{1}} 
\end{equation}
for $t$ large enough. \\
By using the equations $(\ref{eq21})$, $(\ref{eq22})$, $(\ref{eq23})$ and $(\ref{eq24})$, one gets that : 
\begin{equation}
\label{eq25}
dist(Mn,t \mathbb{S}^{1}) \leqslant (\lVert e_{1}(L) \rVert + \lVert e_{2}(L) \rVert) \frac{t}{\lVert L \rVert_{1}} \lVert n \rVert \textit{.}
\end{equation}
\end{proof}
Now we can prove Lemma $\ref{lemme32}$.
\begin{proof}[Proof of Lemma $\ref{lemme32}$]
By using the equation ($\ref{eq19}$), we have that : 
$$\Delta_{1}(L,t) \leqslant \frac{1}{\sqrt{t}}\sum_{n \in \mathbb{Z}^{2}} |\chi_{\mathbb{S}^{1}, M}(n;t) - \mathbf{1}_{t \Omega_{M^{-1}\mathbb{S}^{1}}}(n)| \textit{.}$$
So, Sublemma $\ref{lemme3}$ and Sublemma $\ref{lemme4}$ give that : 
$$\Delta_{1}(L,t) \leqslant \frac{1}{\sqrt{t}} \sum_{n \in \mathbb{Z}^{2}-\{ 0 \}} e^{-C^{2} \frac{t^{4}}{4} \lVert n \rVert^{2} } + \frac{e^{-\frac{t^{4}}{4}}}{\sqrt{t}} = O( \frac{1}{\sqrt{t}})\textit{.} $$
\end{proof}
\subsubsection{Proof of Lemma $\ref{lemme33}$}
To prove Lemma $\ref{lemme33}$, we first need to give another expression of $F(M,t)$, obtained via the Poisson formula. It is the object of the following lemma : 
\begin{lemma}
\label{lemme34}
$$F(M,t) = \frac{1}{\pi}\sum_{n \in \mathbb{Z}^{2}-\{ 0 \}} \frac{1}{\lVert (M^{-1})^{T}n \rVert^{\frac{3}{2}}} \cos( 2 \pi t \lVert (M^{-1})^{T}n \rVert - \frac{3 \pi}{4})  \widetilde{\lambda_{M}}(2 \pi n ; t) + O_{M}(t^{-1}) $$
where the $M$ in index of $O_{M}$ is to signal that it depends on $M$ (or, equivalently, on the lattice $L$). 
\end{lemma}
To prove it, we first need a calculatory sublemma : 
\begin{sublemma}
\label{lemme2}
Let $D \in SL_{2}(\mathbb{R})$ and $\gamma = D \mathbb{S}^{1}$. Then one has for every $\xi \in \mathbb{R}^{2}-\{ 0 \}$ :
$$x_{\gamma}(\xi) = D \frac{D^{T} \xi}{\lVert D^{T} \xi \rVert} \textit{, }$$
$$\rho_{\gamma}(\xi) = \frac{\lVert  \xi \rVert^{3}}{\lVert D^{T} \xi \rVert^{3}} \textit{ and } $$
$$ Y_{\gamma}(\xi) = \lVert D^{T} \xi \rVert $$ 
where $D^{T}$ is the transpose of the matrix $D$.
\end{sublemma}
\begin{proof}
Let's set $ y = D^{-1} x_{\gamma}(\xi) \in \mathbb{S}^{1}$. One knows that : 
\begin{equation}
\label{eq4}
k=y
\end{equation} 
where $k$ is the outer normal to $\mathbb{S}^{1}$ at $y$. \\
Let's call $T$ the unit tangent vector of $\mathbb{S}^{1}$ such that $(k,T)$ is an orthonormal and direct basis of $\mathbb{R}^{2}$. Let's call $R$ the rotation matrix of angle $\frac{\pi}{2}$, so that $$R = \begin{pmatrix} 0 & 1 \\ -1 & 0 \end{pmatrix}$$ and so one has : 
\begin{equation}
\label{eq3}
R k = T \textit{.}
\end{equation}
By definition of $ x_{\gamma}(\xi)$ and because $D$ is a linear map such that $\text{det}(D) > 0$, one knows that $D T$ is a tangent vector of $\gamma$ at the point $x_{\gamma}(\xi)$ such that $(\frac{\xi}{\lVert \xi \rVert}, D T)$ is an orthogonal and direct basis of $\mathbb{R}^{2}$. By using ($\ref{eq3}$), we now know that : $( \xi , D R k)$ is an orthogonal and direct basis of $\mathbb{R}^{2}$.\\
By using the fundamental property of the adjoint of an operator, we get that : 
\begin{equation} 
\label{eq6} (R^{-1} D^{T} \xi , k) \text{ is an orthogonal basis of } \mathbb{R}^{2} \textit{.}
\end{equation}
So one gets that : $R^{-1} D^{T} \xi = a T $ where $a$ is a non-zero real. \\
Finally, we have  
\begin{equation}
\label{eq5}
D^{T} \xi = -a k 
\end{equation}
because of ($\ref{eq3}$). \\
The equation $(\ref{eq4})$ and the fact that $k$ is unitary give us that : 
$$ x_{\gamma}(\xi) = \pm D \frac{D^{T} \xi}{\lVert D^{T} \xi \rVert} \textit{.}$$
Yet, because $( \xi , D R k)$ is a direct basis of $\mathbb{R}^{2}$ we know that : 
$$( \xi , \pm D R D^{T} \xi) \textit{ is a direct basis.}$$
Yet, one can see that $$R D R = - (D^{-1})^{T}$$ and so one gets that : 
$$(R \xi, \mp \xi) \textit{ is a direct basis.} $$
So, one gets finally that : 
\begin{equation}
\label{eq7}
x_{\gamma}(\xi) =  D \frac{D^{T} \xi}{\lVert D^{T} \xi \rVert} \textit{.}
\end{equation}
As a consequence, we get immediately that : 
$$ Y_{\gamma}(\xi) = \lVert D^{T} \xi \rVert \textit{.}$$
Then by using the parametrization $t \longmapsto (\cos(t), \sin(t))$ of the circle and the expression of a curvature radius when using a parametrization, one gets that at a point $X$ of $\gamma = D \mathbb{S}^{1}$, $$\rho(X)=\lVert (D^{-1})^{T} D^{-1} X \rVert^{3} \textit{.} $$
So, with equation $(\ref{eq7})$, finally one gets that : 
$$\rho_{\gamma}(\xi) = \frac{\lVert  \xi \rVert^{3}}{\lVert D^{T} \xi \rVert^{3}} \textit{.} $$
\end{proof}
We can now tackle the proof of Lemma $\ref{lemme34}$.
\begin{proof}[Proof of Lemma $\ref{lemme34}$]
According to the equation ($\ref{eq16}$), the Poisson summation formula and because of the fact that $\tilde{\mathbf{1}_{t \Omega_{M^{-1}\mathbb{S}^{1}}}} (0) = \text{Area}(t \Omega_{M^{-1}\mathbb{S}^{1}})$ one has that : 
\begin{equation}
\label{eq26}
F(M,t) = \frac{1}{\sqrt{t}} \sum_{n \in \mathbb{Z}^{2}-\{ 0 \}} \widetilde{\mathbf{1}_{t \Omega_{M^{-1}\mathbb{S}^{1}}}}(2 \pi n) \widetilde{\lambda_{M}}(2 \pi n ; t) \textit{.}
\end{equation}
Yet, according to Lemma 2.1 from $\cite{bleher1992distribution}$, one has that :
\begin{equation}
\label{eq27}
  \widetilde{\mathbf{1}_{t \Omega_{M^{-1}\mathbb{S}^{1}}}}(\xi) = \sqrt{t} \lVert \xi \rVert^{-\frac{3}{2}} \sum_{\pm} \sqrt{2 \pi \rho_{M^{-1} \mathbb{S}^{1}}( \pm \xi)} \exp (\pm i( t Y_{ M^{-1} \mathbb{S}^{1}}(\pm \xi) - \frac{3 \pi}{4})) + O_{M}(t^{- \frac{1}{2}} \lVert \xi \rVert ^{-\frac{5}{2}})
\end{equation}
By using Sublemma $\ref{lemme2}$, we get with equation ($\ref{eq26}$) : 
\begin{equation}
\label{eq28}
F(M,t) = \frac{1}{\pi}\sum_{n \in \mathbb{Z}^{2}-\{ 0 \}} \frac{1}{\lVert (M^{-1})^{T}n \rVert^{\frac{3}{2}}} \cos( 2 \pi t \lVert (M^{-1})^{T}n \rVert - \frac{3 \pi}{4})   \widetilde{\lambda_{M}}(2 \pi n ; t) + O_{M}(t^{-1}) \textit{.}
\end{equation}
\end{proof}
Using the equation ($\ref{eq13}$), the equation ($\ref{eq19}$), the fact that if $M$ represents a lattice $L$, $(M^{-1})^{T}$ represents the dual lattice $L^{\perp}$, and the previous lemma, that is Lemma $\ref{lemme34}$, one gets that : 
\begin{equation}
\label{eq29}
(\Delta_{2})_{A}(L,t) \leqslant \Delta_{2,1}(L,t) + (\Delta_{2,2})_{A}(L,t) + (\Delta_{2,3})_{A}(L,t)
\end{equation}
where 
\begin{equation}
\label{eq30}
\Delta_{2,1}(L,t) = O_{M}(t^{-1})
\end{equation}
\begin{equation}
\label{eq31}
(\Delta_{2,2})_{A}(L,t) = | \sum_{\substack{ l \in L^{\perp} \\ 0 < \lVert l \rVert \leqslant A }} \frac{1}{\lVert l \rVert^{\frac{3}{2}}} \cos( 2 \pi t \lVert l \rVert - \frac{3 \pi}{4})(1 - e^{- (2 \pi)^{2} \frac{\lVert l \rVert^{2}}{t^{2}}}) |
\end{equation}
\begin{equation}
\label{eq32}
(\Delta_{2,3})_{A}(L,t) =  |\sum_{\substack{ l \in L^{\perp} \\ A < \lVert l \rVert }} \frac{1}{\lVert l \rVert^{\frac{3}{2}}} \cos( 2 \pi t \lVert l \rVert - \frac{3 \pi}{4}) e^{- (2 \pi)^{2} \frac{\lVert l \rVert^{2}}{t^{2}}}| \textit{.}
\end{equation}
So, if we prove the following lemmas, we will get Lemma $\ref{lemme33}$ and, $\textit{in fine}$, get Lemma $ \ref{lemme30}$ : 
\begin{lemma} 
\label{lemme35}
$\Delta_{2,1}(L,t)$ converges almost surely to $0$ when $t \rightarrow \infty$.
\end{lemma}
Let's remark, by the way, that this last lemma is immediate according to equation $(\ref{eq30})$.
\begin{lemma}
\label{lemme36}
For all $A > 0$, $(\Delta_{2,2})_{A}(L,t)$ converges to $0$ when $t \rightarrow \infty$.
\end{lemma}
\begin{lemma}
\label{lemme37}
For all $\alpha > 0$, for all $A$ large enough, for all $t$ large enough,
$$\mathbb{P}((\Delta_{2,3})_{A}(L,t) \geqslant \alpha) \leqslant \alpha \textit{.}$$
\end{lemma}
\begin{proof}[Proof of Lemma $\ref{lemme33}$]
Let $\alpha > 0$. Let's take $A$ large enough so that for all $t$ large enough, $$\mathbb{P}((\Delta_{2,3})_{A}(L,t) \geqslant \alpha) \leqslant \alpha \textit{.}$$ 
It is possible according to Lemma $\ref{lemme37}$. \\
According to Lemma $\ref{lemme36}$, according to Lemma $\ref{lemme35}$ and because the almost-sure convergence imply the convergence in probability, even if it means taking $t$ larger, one can suppose that : $$\mathbb{P}((\Delta_{2,1})(L,t) \geqslant \alpha) \leqslant \alpha \textit{ and}$$ 
$$\mathbb{P}((\Delta_{2,2})_{A}(L,t) \geqslant \alpha) \leqslant \alpha \textit{.}$$ 
By using equation $\ref{eq29}$, one gets the wanted result.
\end{proof}
Before following with the proof of Lemma $\ref{lemme30}$, let's say a few words about Lemma $\ref{lemme36}$ and Lemma $\ref{lemme37}$. The first tells us that the non-regularized Fourier series is "close" enough to the regularized Fourier series whereas the second one tells us that the large terms of the regularized Fourier series do not matter, in a certain sense, for our study. \\
It remains only to prove Lemma $\ref{lemme36}$ and Lemma $\ref{lemme37}$. We will do just that in the next two subsubsections.
\subsubsection{Proof of Lemma $\ref{lemme36}$}
\begin{proof}[Proof of Lemma $\ref{lemme36}$]
Let $l \in L^{\perp}$. Then one has : 
\begin{equation}
\label{eq33}
|1 - e^{- (2 \pi)^{2} \frac{\lVert l \rVert^{2}}{t^{2}}} | \leqslant \frac{(2 \pi)^{2} \lVert l \rVert^{2}}{t^{2}} \textit{.}
\end{equation}
With this equation and with equation $(\ref{eq31})$, one gets that : 
\begin{equation}
\label{eq34}
(\Delta_{2,2})_{A}(L,t) \leqslant \sum_{\substack{l \in L^{\perp} \\ 0 < \lVert l \rVert \leqslant A}} \frac{(2 \pi)^{2}}{t^{2}} \lVert l \rVert^{\frac{1}{2}} \textit{.}
\end{equation}
It follows that there exists $C(L) > 0$ such that : 
\begin{equation}
\label{eq35}
(\Delta_{2,2})_{A}(L,t) \leqslant C(L) \frac{A^{\frac{5}{2}}}{t^{2}} \textit{.}
\end{equation}
\end{proof}
\subsubsection{Proof of Lemma $\ref{lemme37}$} 
To prove Lemma $\ref{lemme37}$ we need to use what are called $\textit{Siegel}$ and $\textit{Rogers formulas}$. Theses formulas will also be useful later in this paper. \\
By setting $c_{k} = \zeta(2)^{-k}$ for $k$ an integer larger than $1$ and where $\zeta$ denotes the $\zeta$ function of Riemann, one has the following formulas :
\begin{lemma}[\cite{marklof1998n},\cite{vinogradov2010limiting},\cite{kelmer2021second}]
\label{lemme5} 
For $f$ a piecewise smooth function with compact support on $\mathbb{R}^{2}$, one has : \begin{itemize}
\item $$  \int_{\mathscr{S}_{2}} \mathcal{S}(f) d \mu_{2} = c_{1} \int_{\mathbb{R}^{2}} f d\lambda  $$
\item When $f$ is even, 
\begin{align*}
(b) \int_{\mathscr{S}_{2}} \mathcal{S}(f)^{2} d \mu_{2}  \leqslant C \int_{\mathbb{R}^{2}} f^{2} d\lambda + c_{2}(\int_{\mathbb{R}^{2}}f d\lambda)^{2} 
\end{align*}
where $C > 0$. 
\end{itemize}

\end{lemma}
With this lemma, we are going to prove two lemmas that will enable us to prove Lemma $\ref{lemme37}$ by using Chebyshev's inequality : the first one is intended to estimate the expectation of $(\Delta_{2,3})_{A}(L,t)$ to see that it goes to $0$ when $t \rightarrow \infty$ (uniformly in $A$), the second one is intended to estimate its variance to see that it can be as uniformly small in $t$ as one wants if $A$ is chosen large enough. Until the end of this section, we are going to suppose $A > 1$.
\begin{lemma}
\label{lemme6}
$$\mathbb{E}((\Delta_{2,3})_{A}(L,t)) = O(\frac{1}{t }) \textit{.}$$
\end{lemma}
\begin{proof}
One has :
\begin{equation}
\label{eq36}
(\Delta_{2,3})_{A}(L,t) =  \sum_{\substack{ l \in L^{\perp} \textit{ prime } \\ A < \lVert l \rVert }} f(l) 
\end{equation}
where 
\begin{equation}
\label{eq37}
f(l) = \frac{1}{\lVert l \rVert^{\frac{3}{2}}} \sum_{k \in \mathbb{N}-\{ 0 \}} \frac{\cos( 2 k \pi t \lVert l \rVert - \frac{3 \pi}{4})}{k^{\frac{3}{2}}} e^{- (2 k \pi)^{2} \frac{\lVert l \rVert^{2}}{t^{2}}}\textit{.}
\end{equation}
The Lemma $\ref{lemme5}$ gives us then that : 
\begin{equation}
\label{eq38}
\mathbb{E}((\Delta_{2,3})_{A}(L,t)) \leqslant C \int_{\mathbb{R}^{2}} f(x) \mathbf{1}_{ \lVert x \rVert > A } dx 
\end{equation}
where $C > 0$ because the density of $\tilde{\mu}_{2}$ is supposed to be bounded. \\
By passing into polar coordinates, one gets that : 
\begin{equation}
\label{eq39}
\mathbb{E}((\Delta_{2,3})_{A}(L,t)) \leqslant C \int_{r=A}^{\infty} \frac{\sum_{k \in \mathbb{N}-\{ 0 \}} \frac{\cos( 2 k \pi t r - \frac{3 \pi}{4})}{k^{\frac{3}{2}}} e^{- (2 k \pi)^{2} \frac{r^{2}}{t^{2}}}}{r^{\frac{1}{2}}}  dr
\end{equation}
(the constant $C$ has changed but it does not matter). 
Lebesgue's dominated convergence theorem gives us that : 
\begin{equation}
\label{eq40}
\mathbb{E}((\Delta_{2,3})_{A}(L,t)) \leqslant C \sum_{k \in \mathbb{N}-\{ 0 \}} \frac{1}{k^{\frac{3}{2}}} \int_{r=A}^{\infty}  \frac{\cos( 2 k \pi t r - \frac{3 \pi}{4})}{r^{\frac{1}{2}}} e^{- (2 k \pi)^{2} \frac{r^{2}}{t^{2}}}dr\textit{.}
\end{equation}
Finally, an integration by part gives us that : 
\begin{align}
\label{eq41}
\mathbb{E}((\Delta_{2,3})_{A}(L,t))\leqslant C \sum_{k \in \mathbb{N}-\{ 0 \}} \frac{1}{k^{\frac{3}{2}}}( & - \frac{\sin(2 k \pi t A - \frac{3 \pi}{4}) e^{-(2 k \pi)^{2} \frac{A^{2}}{t^{2}}}}{2 k \pi t A^{\frac{1}{2}}} \nonumber \\
 & + \int_{A}^{\infty} \frac{\sin(2 k \pi t A - \frac{3 \pi}{4}) e^{-(2 k \pi)^{2} \frac{A^{2}}{t^{2}}}}{4 k \pi t r^{\frac{3}{2}}} dr \nonumber \\
 & + \int_{A}^{\infty} \frac{\sin (2 k \pi t r - \frac{3 \pi}{4})}{2 k \pi t} \frac{(2 \pi k)^{2}}{t^{2}} 2 r^{\frac{1}{2}} e^{-(2 k \pi)^{2} \frac{r^{2}}{t^{2}}} dr) \textit{.}
\end{align}
By using that $r^{\frac{1}{2}} \leqslant r$ (because $A > 1$), one has finally, by estimating the three terms of the right member : 
$$\mathbb{E}((\Delta_{2,3})_{A}(L,t))= O(\frac{1}{t}) \textit{.}$$
\end{proof}
\begin{lemma}
\label{lemme7}
$$\text{Var}((\Delta_{2,3})_{A}(L,t)) = O(\frac{1}{A})  $$
where the $O$ can be chosen independent from $t$. 
\end{lemma}
\begin{proof}
By using the same notation as before, by using again the Lemma $\ref{lemme5}$ and by using the Lemma $\ref{lemme6}$, one gets that : 
\begin{equation}
\label{eq42}
\text{Var}((\Delta_{2,3})_{A}(L,t)) \leqslant C \int_{\mathbb{R}^{2}} f^{2}(x) \mathbf{1}_{\lVert x \rVert > A} dx  \textit{.}
\end{equation}
So, by passing into polar coordinates, one gets that : 
\begin{equation}
\label{eq43}
\text{Var}((\Delta_{2,3})_{A}(L,t))\leqslant C 2 \pi \int_{r=A}^{\infty} \frac{1}{r^{2}} (\sum_{k \in \mathbb{N}-\{ 0 \}} \frac{\cos( 2 k \pi t r - \frac{3 \pi}{4})}{k^{\frac{3}{2}}} e^{- (2 k \pi)^{2} \frac{r^{2}}{t^{2}}})^{2} dr  \textit{.}
\end{equation}
Because $\sum_{k=1}^{\infty} \frac{1}{k^{\frac{3}{2}}} < \infty$, we get the wanted result. 
\end{proof}
We can now prove the Lemma $\ref{lemme37}$.
\begin{proof}[Proof of Lemma $\ref{lemme37}$]
The Chebyshev's inequality gives the wanted result if, first, we choose $A$ large enough and, second, we choose $t$ large enough so that $\mathbb{E}((\Delta_{2,3})_{A}(L,t))$ and $\text{Var}((\Delta_{2,3})_{A}(L,t))$ are small enough. These choices are possible according to Lemmas $\ref{lemme6}$ and $\ref{lemme7}$.
\end{proof}
So, now the proof of Lemma $\ref{lemme30}$ is complete and we will conclude this section by proving the Lemma $\ref{lemme31}$ so that the proof of Proposition $\ref{prop8}$ will be complete.
\subsection{Proof of Lemma $\ref{lemme31}$}
To prove the Lemma $\ref{lemme31}$, we are going to take the same kind of approach as before : estimate the expectation and the variance of the quantity $S_{A}-S_{A,prime}$ and get the result via the Chebyshev's inequality. \\
We have that : 
\begin{equation}
\label{eq45}
H_{A}(L,t)-S_{A,prime}(L,t) = \sum_{l \in L^{\perp} \textit{prime}} f(l) 
\end{equation}
where
\begin{equation}
\label{eq46}
f(l) = \frac{1}{ \lVert l \rVert^{\frac{3}{2}}} \mathbf{1}_{0 < \lVert l \rVert \leqslant A} \sum_{k \geqslant \lfloor \frac{A}{\lVert l \rVert} \rfloor + 1} \frac{\cos(2 k \pi t \lVert l \rVert - \frac{3 \pi}{4})}{k^{\frac{3}{2}}} \textit{.}
\end{equation}
With this expression, we see that we are going to have a little problem of integrability at $0$ when using the Lemma $\ref{lemme5}$. That's why, we have to exclude $0$ and we will suppose that $L$ is chosen so that $\lVert L^{\perp} \rVert_{1} \geqslant \epsilon$ where $0 < \epsilon < 1$. Only a small number of lattices is excluded according to this lemma that we recall : 
\begin{lemma}
\label{lemme8}
For every $0 < \epsilon < 1$, one has that 
$$\mathbb{P}(\lVert L \rVert_{1} < \epsilon) = O(\epsilon^{2}) \textit{.} $$
\end{lemma}
\begin{proof}
It is a consequence of Lemma $\ref{lemme5}$ by taking 
$$\mathcal{S}(f)(L) = \sum_{l \in L} \mathbf{1}_{B_{f}(0,\epsilon)}(l) $$ 
where $ \mathbf{1}_{B_{f}(0,\epsilon)}(l)$ is the indicator function of the closed ball for the norm $\lVert \rVert$ centred on $0$ of radius $\epsilon$.
\end{proof}
Thus, for the chosen lattices, we have :
\begin{equation}
\label{eq47}
H_{A}(L,t)-S_{A,prime}(L,t) = \sum_{l \in L^{\perp} \textit{prime}} f(l) \mathbf{1}_{\lVert l \rVert \geqslant \epsilon} = \Delta_{3,\epsilon,A,t}(L) \textit{.}
\end{equation}
\begin{lemma}
\label{lemme9}
$$\mathbb{E}(\Delta_{3,\epsilon,A,t}(L)) = O_{\epsilon,A}(\frac{1}{t}) \textit{.}$$
\end{lemma}
\begin{proof}
By using the Lemma $\ref{lemme5}$, one gets that : 
\begin{equation}
\label{eq48}
\mathbb{E}(\Delta_{3,\epsilon,A,t}(L)) \leqslant C \int_{\epsilon}^{A} \frac{1}{r^{\frac{1}{2}}} \sum_{k \geqslant \frac{A}{r}} \frac{\cos(2 k \pi t r - \frac{3 \pi}{4})}{k^{\frac{3}{2}}} dr \textit{.} 
\end{equation}
Lebesgue's dominated convergence theorem gives us that : 
\begin{equation}
\label{eq49}
\int_{\epsilon}^{A} \frac{1}{r^{\frac{1}{2}}} \sum_{k \geqslant \frac{A}{r}} \frac{\cos(2 k \pi t r - \frac{3 \pi}{4})}{k^{\frac{3}{2}}} dr = \sum_{k=1}^{\infty} \frac{1}{k^{\frac{3}{2}}} \int_{\max(\frac{A}{k}, \epsilon)}^{A} \frac{\cos(2 k \pi t r - \frac{3 \pi}{4})}{r^{\frac{1}{2}}}dr \textit{.}
\end{equation}
An integration by part as in the proof of lemma $\ref{lemme6}$ and the equation $(\ref{eq48})$ give us finally that :  
\begin{equation}
\label{eq50}
\mathbb{E}(\Delta_{3,\epsilon,A,t}(L)) = O_{\epsilon,A}(\frac{1}{t}) \textit{.}
\end{equation}
\end{proof}
\begin{lemma}
\label{lemme10}
There exists $K > 0$ such that : 
$$\text{Var}(\Delta_{3,\epsilon,A,t}(L)) \leqslant K(- \frac{\log(\epsilon)}{A} + \frac{\log(A)}{A}) \textit{.} $$
\end{lemma}
\begin{proof}
Lemma $\ref{lemme5}$ gives us that : 
\begin{equation}
\label{eq51}
\text{Var}(\Delta_{3,\epsilon,A,t}(L)) \leqslant C 2 \pi \int_{\epsilon}^{A} \frac{1}{r^{2}} (\sum_{k \geqslant \frac{A}{r}} \frac{\cos(2 \pi k t r - \frac{3 \pi }{4})}{k^{\frac{3}{2}}})^{2} dr \textit{.}
\end{equation}
Yet one also has that for all $x > 0$ : 
\begin{equation}
\label{eq52}
\sum_{k \geqslant x} \frac{1}{k^{\frac{3}{2}}} \leqslant \frac{D}{x^{\frac{1}{2}}}
\end{equation}
where $D >0$. \\
Thus, equations $(\ref{eq51})$ and $(\ref{eq52})$ imply that : 
\begin{equation}
\label{eq53}
\text{Var}(\Delta_{3,\epsilon,A,t}(L)) \leqslant \frac{2 \pi C D}{A} \int_{\epsilon}^{A} \frac{1}{r} dr \textit{.}
\end{equation}
\end{proof}
We can now give the proof of Lemma $\ref{lemme31}$.
\begin{proof}[Proof of Lemma $\ref{lemme31}$]
First we take $1> \epsilon > 0$ small enough so that the measure of the neglected lattices, $\textit{id est}$ the lattices such that $\lVert L \rVert_{1} < \epsilon$, is small enough. It is possible according to Lemma $\ref{lemme8}$. \\
Then we take $A$ large enough so that $\text{Var}(\Delta_{3,\epsilon,A,t}(L))$ is small enough. It is possible according to Lemma $\ref{lemme10}$. \\
Finally, we take $t$ large enough so that $\mathbb{E}(\Delta_{3,\epsilon,A,t}(L))$ is small enough, which is possible according to Lemma $\ref{lemme9}$, and conclude by using Chebyshev's inequality.
\end{proof}
So, we are now brought back to the study of $S_{A,prime}(L,t)$ when $t \rightarrow \infty$ and the next section is dedicated to it. \\
$\textbf{We are going to replace}$ $L^{\perp}$ $\textbf{by}$ $L$ (it changes nothing because we are studying the asymptotic convergence in law with $L \in \mathscr{S}_{2}$ distributed according to $\tilde{\mu}_{2}$).
\section{Study of $S_{A,prime}(L,t)$ when $t \rightarrow \infty$}
\subsection{Reductions for the study of $S_{A,prime}(L,t)$ and proof of Theorem $\ref{thm2}$}
Before entering in the main object of this section, we need to do a small rewriting of $S_{A,prime}(L,t)$.\\
We recall that a vector $l \in L$ is prime if, and only if, $K l \in K L$ is prime where $K \in SL_{2}(\mathbb{R})$. Furthermore, a vector $(l_{1},l_{2}) \in \mathbb{Z}^{2}$ is prime if, and only if, $l_{1} \wedge l_{2} = 1$. \\
By using the symmetry $l \longmapsto -l$, we deduce that $S_{A,prime}(L,t)$ can be rewritten as followed : 
\begin{equation}
\label{eq55}
S_{A,prime}(L,t) = \frac{2}{\pi} \sum_{k \in \Pi_{A}(L)} \frac{Z_{k}(L,t)}{Y_{k}(L)} 
\end{equation} 
where, for $k=(k_{1},k_{2}) \in \mathbb{Z}^{2}$,
\begin{equation}
\label{eq56}
Y_{k}(L) = \lVert k_{1} e_{1}(L) + k_{2} e_{2}(L) \rVert^{\frac{3}{2}} \textit{, }
\end{equation}
\begin{equation}
\label{eq57}
Z_{k}(L,t) = \sum_{m \in \mathbb{N}-\{ 0 \}} \frac{ \cos(2 \pi t m \lVert k_{1} e_{1}(L) + k_{2} e_{2}(L) \rVert - \frac{3 \pi}{4})}{m^{\frac{3}{2}}} \textit{, }
\end{equation}
and where
\begin{equation}
\label{eq58}
\Pi_{A}(L)= \{ (k_{1},k_{2}) \in \mathbb{Z}^{2} \text{ } | \text{ } k_{1} \wedge k_{2} = 1 \textit{, } k_{1} \geqslant 0 \textit{ } \lVert k_{1} e_{1}(L) + k_{2} e_{2} (L) \rVert \leqslant A \} 
\end{equation}
and here we agree that if $(k_{1},k_{2}) \in \Pi_{A}(L)$ then $k_{1} = 0$ implies that $k_{2} =1$ (for the definition of $e_{1}(L)$ and $e_{2}(L)$ see Definition $\ref{def1}$). \\
\\
Let's recall that : 
\begin{equation}
\label{eq67}
\Pi = \{ (k_{1},k_{2}) \in \mathbb{Z}^{2} \text{ } | \text{ } k_{1} \wedge k_{2} = 1 \textit{, } k_{1} \geqslant 0  \} \textit{.}
\end{equation}
Our goal now is to prove the following proposition : 
\begin{prop}
\label{prop10}
$ \{  Z_{k}(L,t) \}_{k \in \Pi}$ converge, when $t \rightarrow \infty$, towards independent identically distributed real random variables that have a compact support, are symmetrical and are independent of $L$. 
\end{prop}
In the next section we are going to consider the sums of the type $$\tilde{S}_{A}(\omega,L) = \sum_{k \in \Pi_{A}(L)} \frac{Z_{k}(\omega)}{\lVert k_{1} e_{1}(L) + k_{2} e_{2}(L) \rVert^{\frac{3}{2}}}   $$
where $Z_{k}$ are non-zero real independent identically distributed random variables from $\Omega \ni \omega$ that are symmetrical and have a compact support. We will see, with Proposition $\ref{prop15}$, that the sums of this type : \begin{itemize}
\item converge almost surely (on the space $\Omega \times \mathscr{S}_{2}$)
\item their respective limits are symmetrical and their expectations are equal to $0$
\item their respective limits admit moment of order $1+\kappa$ for any $0 \leqslant \kappa < \frac{1}{3}$ 
\item their respective limits do not admit a moment of order $\frac{4}{3}$ when $\sigma(L) \geqslant m$ where $m > 0$ and where $L$ belongs to an event of the form $( \lVert L \rVert < \alpha )$ with $\alpha > 0$. 
\item when there exists $\alpha > 0 $ such that $\tilde{\mu}_{2}(\{ L \in \mathscr{S}_{2} \text{ } | \text{ } \lVert L \rVert_{1} < \alpha \}) = 0$ then their respective limits admit moments of all order $1 \leqslant p < \infty$.  
\end{itemize}
We are going to see now that it is enough to prove Proposition $\ref{prop10}$ and Proposition $\ref{prop15}$ to establish Theorem $\ref{thm2}$, with the exception of the exact form of the limiting law (yet it is given by Proposition $\ref{prop11}$).
\begin{proof}[Proof of Theorem $\ref{thm2}$]
Let $\psi \in C_{c}^{\infty}(\mathbb{R})$. Let $\epsilon > 0$. According to Proposition $\ref{prop8}$, we can take $A$ as large as we want and then $t$ as large as we want so that : 
\begin{equation}
\label{eq502}
| \mathbb{E}\left( \psi(\frac{\mathcal{R}(t \mathbb{D}^{2}, L)}{\sqrt{t}}) \right) - \mathbb{E}\left( \psi(S_{A,prime}(L,t)) \right) | \leqslant \epsilon \textit{.} 
\end{equation}
Furthermore, even if it means neglecting lattices that form a set of small measure, one can suppose that only a finite number of $k$, independent of $L$, intervene in the equation $(\ref{eq55})$. \\
Thus, thanks to $\ref{prop10}$, one has that : 
\begin{equation}
\label{eq503}
| \mathbb{E}\left( \psi(S_{A,prime}(L,t)) \right) - \mathbb{E}\left( \psi(\tilde{S}_{A}(\omega,L)) \right) | \leqslant \epsilon 
\end{equation}
where the $Z_{k}(\omega)$ in $\tilde{S}_{A}(\omega,L)$ are given by Proposition $\ref{prop10}$. \\
Furthermore, Proposition $\ref{prop15}$ gives us that : 
\begin{equation}
\label{eq504}
| \mathbb{E}\left( \psi(\tilde{S}_{A}(\omega,L)) \right) - \mathbb{E}\left( \psi(\lim_{A \rightarrow \infty} \tilde{S}_{A}(\omega,L)) \right) | \leqslant \epsilon 
\end{equation}
with $\lim_{A \rightarrow \infty} \tilde{S}_{A}(\omega,L))$ that verify all the listed properties. \\
So, equations $(\ref{eq502})$, $(\ref{eq503})$ and $(\ref{eq504})$ give the wanted result.
\end{proof}
The main reason why $Z_{k}$ are going to be independent from $L$ is the presence of the factor $t$. \\
The main reasons why the rest of Proposition $\ref{prop10}$ will be true are the presence of the factor $t$ in $Z_{k}$ and the fact that the coefficients of order 1 of the Taylor series of $(\lVert k_{1} e_{1}(L) + k_{2} e_{2}(L) \rVert)_{k \in \Pi_{A}(L)}$ on a small geodesic segment are $\mathbb{Z}$-free. \\
\\
In order to prove Proposition $\ref{prop10}$, it is actually enough to prove the following proposition : 
\begin{prop}
\label{prop11}
For $k=(k_{1},k_{2})$, let 
\begin{equation}
\label{eq59}
\theta_{k}(L,t) = t  \lVert k_{1} e_{1}(L) + k_{2} e_{2}(L) \rVert \mod 1 \textit{.}
\end{equation}
Then, we have that
$\{ \theta_{k}(L,t) \}_{k \in \Pi}$ converge, when $t \rightarrow \infty$, towards random variable that are independent identically distributed, are distributed according to the Lebesgue measure $\lambda$ over $\mathbb{R}/\mathbb{Z}$ and are independent from $L$.
\end{prop}
$\textbf{Thanks to this proposition, we now understand why the limit law of }$ $ \frac{\mathcal{R}(t \mathbb{D}^{2},L)}{\sqrt{t}}$ $ \textbf{ is} $ \\
$\textbf{given by } S(\theta,L) \textit{.}$
To prove this last proposition, it is sufficient to prove the following proposition where $e(\theta)$ stands for $\exp(i 2 \pi \theta)$ : 
\begin{prop}
\label{prop12}
For every $l \in \mathbb{N}-\{0\}$, for every $\psi \in C_{c}^{\infty}(\mathscr{S}_{2})$, for every $(p_{1},\cdots,p_{l}) \in \mathbb{Z}^{l}-\{ 0 \}$, one has : 
\begin{equation}
\label{eq66}
\mathbb{E} \left( \psi(L) e(\sum_{h=1}^{l} p_{h} \theta_{k_{h}}) \right) \underset{t \rightarrow \infty}{\rightarrow} 0
\end{equation}
where the $k_{h} \in \Pi$ are all distinct. 
\end{prop}
Before passing to the proof of Proposition $\ref{prop12}$, let's give some heuristic about it.\\ Basically, by working with a foliation of the space $\mathscr{S}_{2}$ given by small enough geodesic segments, we are first going to have : $$\mathbb{E} \left(  \psi(L) e(\sum_{h=1}^{l} p_{h} \theta_{k_{h}}) \right) \approx \mathbb{E}(\psi) \mathbb{E} \left(e(\sum_{h=1}^{l} p_{h} \theta_{k_{h}}) \right)  $$
due to the presence of the factor $t$ in $\theta$. \\
The right member will go to $0$ when $t$ goes to infinity because a Riemann-Lebesgue lemma will apply because quantities "close" to the variables $\theta_{k}$ are typically $\mathbb{Z}$-free (see the heuristic explanation of the second step). \\ 
The rest of this section is now dedicated to the proof of the Proposition $\ref{prop12}$.
\subsection{Foliation and local estimates}
We recall that a foliation of the space $\mathscr{S}_{2}$ is given by the orbits of the group $\delta$ where $$ \delta(\lambda) = \begin{pmatrix} \lambda & 0 \\ 0 & \frac{1}{\lambda} \end{pmatrix} \textit{.}$$
To prove Proposition $\ref{prop12}$, we are going to look at what it is happening on a small "segment" of the form 
\begin{equation}
\label{eq60}
S_{\epsilon}(L) = \{ \delta(\lambda) L \textit{ } | \textit{ } \lambda \in [\frac{1}{1+ \epsilon},1+\epsilon] \}
\end{equation}
where $L \in \mathscr{S}_{2}$ and $\epsilon > 0$ can be taken as small as possible. More precisely, we are going to show, when $t \rightarrow \infty$, the independence of the $(\theta_{k})$ and of $L$ over smalls segments of the form $S_{\epsilon}(L)$, as well as the fact that the $(\theta_{k})$ are identically distributed and distributed according to the normalized Haar measure over $\mathbb{R}/ \mathbb{Z}$.  \\
Let's set, for $k=(k_{1},k_{2}) \in \Pi$ : 
\begin{equation}
\label{eq400}
W_{k}(L) = \sum_{j=1}^{2}k_{j}^{2} \frac{(e_{j}(L))_{1}^{2}-(e_{j}(L))_{2}^{2}}{\lVert k_{1} e_{1}(L) + k_{2} e_{2}(L) \rVert} \textit{.}
\end{equation}
On a segment of the form $S_{\epsilon}(L)$, the following lemma basically tells us how we can estimate the quantities $Z_{k}$ : 
\begin{lemma}
\label{lemme11}
For a typical $L \in \mathscr{S}_{2}$, there exists $\epsilon > 0$ small enough such that for every $\lambda \in [\frac{1}{1+\epsilon},1+\epsilon]$, $$e_{1}(\delta(\lambda)L) = \delta(\lambda)e_{1}(L) \textit{ and } e_{2}(\delta(\lambda)L) = \delta(\lambda)e_{2}(L) \textit{.}$$
Furthermore, for such a $\lambda$, for $k=(k_{1},k_{2}) \in \Pi$, we have for $h=\lambda-1$,
\begin{align}
\label{eq61}
& \lVert k_{1} e_{1}(\delta(\lambda)L) + k_{2} e_{2}((\delta(\lambda)L) \rVert = \lVert k_{1} e_{1}(L) + k_{2} e_{2}(L) \rVert \\
& + W_{k}(L) h + O_{k_{1},k_{2},L}(h^{2})  \textit{.} \nonumber 
\end{align}

\end{lemma}
\begin{proof}
First, if the first part of the lemma is acquired, then the following equation are acquired $\textit{via}$ a basic calculation of Taylor series. Thus, we just have to prove the first part of the Lemma $\ref{lemme11}$. \\
Now let's prove that for a $\epsilon > 0$ small enough, for every $\lambda \in [\frac{1}{1+ \epsilon}, 1 + \epsilon]$, $e_{1}(\delta(\lambda)L) = \delta_{\lambda}e_{1}(L) \textit{.}$ \\
We note that for a typical lattice $L$, we have that 
\begin{equation}
\label{eq63}
\min_{l \in L-\{ \pm e_{1}(L), \pm e_{2}, 0 \}} \lVert l \rVert \geqslant \lVert e_{2}(L) \rVert  > \lVert e_{1}(L) \rVert \textit{.}
\end{equation}
So there exists $\kappa > 0$ such that for every $k \in L-\{\pm e_{1}(L), \pm e_{2}(L),0 \}$,
\begin{equation}
\label{eq64}
\lVert k \rVert \geqslant \lVert e_{1}(L) \rVert + \kappa \textit{.}
\end{equation}
This last equation and the continuity of $\lVert \rVert$ give that, for a $\epsilon > 0$ small enough, for every $\lambda \in [\frac{1}{1+ \epsilon}, 1 + \epsilon]$ : 
\begin{equation}
\label{eq65}
\lVert \delta(\lambda) k \rVert \geqslant \lVert \delta(\lambda) e_{1}(L) \rVert + \frac{\kappa}{2} \textit{.}
\end{equation}
This last equation and the fact $(\delta(\lambda) e_{1}(L))_{1} > 0$ give us that $\delta(\lambda) e_{1}(L) = e_{1}(\delta(\lambda) L)$. \\
Finally, let's prove that, even if it means reducing $\epsilon > 0$, for every $\lambda \in [\frac{1}{1+ \epsilon}, 1 + \epsilon]$, $e_{2}(\delta(\lambda)L) = \delta_{\lambda}e_{2}(L) \textit{.}$\\
We have the wanted result when the first inequality in the equation ($\ref{eq63}$) is strict (by reasoning the same way as before). \\
So let's suppose that this last inequality is an equality. Let's call $l$ such that $ l_{1} \geqslant 0$ and such that $\lVert e_{2}(L) \rVert = \lVert l \rVert$. If $l_{1} > 0$ it means that the couple $(e_{1}(L),e_{2}(L))$ is not well-defined. If $l_{1} = 0$ then $L$ belongs to a negligible set according to the Lemma 4.5 from the article $\cite{Skriganov}$.  
\end{proof}
To prove Proposition $\ref{prop12}$ we see, in light of Lemma $\ref{lemme11}$, that it would be convenient to prove the following proposition : 
\begin{prop}
\label{prop13}
For a typical $L \in \mathscr{S}_{2}$, for every $m \in \mathbb{N}-\{ 0 \}$, for every family $(p_{1},\cdots,p_{m}) \in \mathbb{Z}^{m}$, for every $k_{1},\cdots,k_{m} \in \Pi$, all distinct if 
\begin{equation}
\label{eq505}
\sum_{i=1}^{m}  p_{i} W_{k_{i}}(L) = 0
\end{equation}
then $p_{1} = \cdots = p_{m} = 0$. \\
In other words, for a typical $L \in \mathscr{S}_{2}$, $$\left( W_{k_{i}}(L) \right)_{i \geqslant 1}$$ is a $\mathbb{Z}$-free family.
\end{prop}
The next subsection is dedicated to prove this proposition. 
\subsection{Proof of Proposition $\ref{prop13}$}
Let us set : 
\begin{equation}
\label{eq68}
X_{1} = \lVert e_{1}(L) \rVert ^{2} \textit{, }
\end{equation}
\begin{equation}
\label{eq69}
X_{2} = < e_{1}(L),e_{2}(L) > \textit{, }
\end{equation}
\begin{equation}
\label{eq70}
X_{3} = \lVert e_{2}(L) \rVert^{2} \textit{ and }
\end{equation}
\begin{equation}
\label{eq71}
X_{4} = \frac{ \cos(\alpha) }{\sin (\alpha) } 
\end{equation}
where $\beta \in ]- \pi, \pi[ -\{ 0 \}$ is the oriented angle from $e_{1}(L)$ to $e_{2}(L)$ and where $\alpha$ is the oriented angle from $(1,0)$ to $e_{1}(L)$ (see the figure just after). So one has typically $\alpha \in ]- \frac{\pi}{2}, \frac{\pi}{2} [ - \{ 0 \}$ and $\beta \in ]- \frac{\pi}{2}- \alpha, \frac{\pi}{2} - \alpha[- \{0 \}$.
\begin{center}
\includegraphics[scale=0.5]{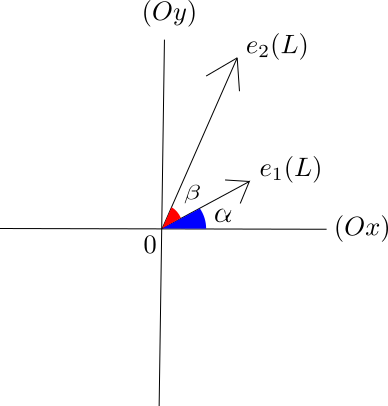}
\end{center}
With theses notations, one has : 
\begin{lemma}
\label{lemme12}
\begin{itemize}
\item $$X_{1} X_{3} = 1 + X_{2}^{2} $$
\item $$ X_{2} = \pm \frac{\cos(\beta)}{\sin(\beta)} $$
\item For every $k \in \Pi$, $$\lVert k_{1} e_{1}(L) + k_{2} e_{2}(L) \rVert = \sqrt{ k_{1}^{2} X_{1} + 2 k_{1} k_{2} X_{2} + k_{2}^{2} X_{3}}$$
\item $$(e_{1}(L))_{1}^{2}-(e_{1}(L))_{2}^{2} = X_{1} \frac{X_{4}^{2} - 1}{X_{4}^{2}+1} $$
\item $$(e_{2}(L))_{1}^{2}-(e_{2}(L))_{2}^{2} =  \frac{1}{X_{1}(1+X_{4}^{2})} ((X_{4}^{2}-1)(X_{2}^{2}-1)-4 \text{ sgn}(\beta) X_{2} X_{4}) $$
\end{itemize}
\end{lemma}
\begin{proof}
Because $L \in \mathscr{S}_{2}$, we have that
\begin{equation}
\label{eq72}
1 = \det(e_{1}(L),e_{2}(L))^{2} = X_{1} X_{3} \sin(\beta)^{2} \textit{.} 
\end{equation}
Furthermore, one has : 
\begin{equation}
\label{eq73}
X_{2}^{2} = X_{1} X_{3} \cos( \beta)^{2} \textit{.} 
\end{equation}
So by summing the equation ($\ref{eq72}$) with the equation ($\ref{eq73}$), one gets that : 

$$
X_{1} X_{3} = 1 + X_{2}^{2} \textit{.} 
$$
The second item is also a direct consequence of equations $(\ref{eq72})$ and $(\ref{eq73})$. \\
\\
The third item is obtained by using the fact that $\lVert \cdot \rVert = \sqrt{< \cdot, \cdot>}$. \\
\\
The fourth item is obtained by using the two following equations : 
\begin{equation}
\label{eq74}
(e_{1}(L))_{1}^{2} = X_{1} \cos^{2}(\alpha) \textit{ and } (e_{1}(L))_{2}^{2} = X_{1} \sin^{2}(\alpha) \textit{.}
\end{equation}
Concerning the fifth and last item, thanks to the two following equations
\begin{equation}
\label{eq75}
(e_{2}(L))_{1}^{2} = X_{3} \cos^{2}(\alpha + \beta) \textit{ and } (e_{2}(L))_{2}^{2} = X_{3} \sin^{2}(\alpha + \beta) 
\end{equation}
we get that :
\begin{equation}
\label{eq76}
(e_{2}(L))_{1}^{2}-(e_{2}(L))_{2}^{2} =  \frac{1}{X_{1}(1+X_{4}^{2})} ((X_{4}^{2}-1)(2X_{2}^{2}-X_{1}X_{3})-4 \text{sgn}(\beta) X_{2} X_{4} \sqrt{X_{1}X_{3}-X_{2}^{2}}) \textit{.}
\end{equation}
We conclude the proof of the fifth item by using the first item that has been proved at the beginning. 
\end{proof}
Now, let's set, for every $k \in \Pi$ 
\begin{equation}
\label{eq77}
A_{k} (X_{1},X_{2})  =  k_{1}^{2} X_{1}^{2} + 2 k_{1} k_{2} X_{1} X_{2} + k_{2}^{2} (1+X_{2}^{2}) \textit{.}
\end{equation}
Thanks to the second equation of Lemma $\ref{lemme12}$, we have that : 
\begin{equation}
\label{eq78}
\lVert k_{1} e_{1}(L) + k_{2} e_{2}(L) \rVert \sqrt{X_{1}}  =  \sqrt{A_{k} (X_{1},X_{2})}
\end{equation}
and so, by multiplying, the equation $(\ref{eq505})$ is equivalent to (with the same notations) : 
\begin{equation}
\label{eq79}
\sum_{i=1}^{m} p_{i} \sqrt{\prod_{h \neq i} A_{k_{h}}}\left( (k_{i})_{1}^{2}  X_{1}^{2} (X_{4}^{2} - 1) +(k_{i})_{2}^{2} ((X_{4}^{2}-1)(X_{2}^{2}-1)-4 \text{ sgn}(\beta) X_{2} X_{4})\right) = 0 \textit{.}
\end{equation}
At the end of the day, according to this last remark, to prove the Proposition $\ref{prop13}$, it is enough to prove the following proposition :
\begin{prop}
\label{prop14}
For a typical $L \in \mathscr{S}_{2}$, for all $A_{k_{1}},\cdots,A_{k_{n}},A_{h_{1}},\cdots,A_{h_{r}}$ where $k_{1},\cdots,k_{n}$ and $h_{1},\cdots,h_{r}$ are all distinct from one another and from $(1,0)$, the square root $\sqrt{\prod_{i=1}^{r} A_{h_{i}}} $ does not belong to the field $\mathbb{Q}(X_{1},X_{2},X_{4})[\sqrt{A_{k_{1}}},\cdots,\sqrt{A_{k_{n}}}]$ with $n \geqslant 0$, $r \geqslant 1$.
Consequently, the family formed by the square roots of product of different $A_{k}$ are $\mathbb{Q}(X_{1},X_{2},X_{4})$-free. 
\end{prop}
We must exclude the case $(1,0)$ because $A_{(1,0)} = X_{1}^{2}$. Yet, in the proof of the Proposition $\ref{prop13}$, there won't be any major difficulties created because of this last fact. \\
This proposition is an analogue of a lemma that can be used to prove that the family $(\sqrt{p})_{p \in \mathcal{P}}$, where $\mathcal{P}$ is the set of prime numbers, is $\mathbb{Z}$-free. \\
The Proposition $\ref{prop14}$ will be proved by induction but to achieve it we first need a lemma.
\begin{lemma}
\label{lemme13}
For a typical lattice $L \in \mathscr{S}_{2}$, for all $P \in \mathbb{Q}[Y_{1},Y_{2},Y_{3}]$, $$\textit{if } P(\lVert e_{1}(L) \rVert^{2}, \frac{\cos(\beta)}{\sin(\beta)}, \frac{\cos(\alpha)}{\sin(\alpha)}) = 0 \textit{ then } P = 0 \textit{.} $$
As a consequence, for a typical lattice $L \in \mathscr{S}_{2}$, for all $P \in \mathbb{Q}[Y_{1},Y_{2},Y_{3}]$, $$\textit{if } P(X_{1}, X_{2}, X_{4}) = 0 \textit{ then } P = 0  $$
where $X_{1}$ was defined by equation $(\ref{eq68})$, $X_{2}$ was defined by equation $(\ref{eq69})$ and $X_{4}$ was defined by equation $(\ref{eq71})$.
\end{lemma}
\begin{proof}
The result is immediate if one has in mind that : 
\begin{itemize}
\item A typical lattice $L \in \mathscr{S}_{2}$ is determined by three, and only three, continuous parameters : $\alpha \in ]- \frac{\pi}{2}, \frac{\pi}{2}[-\{ 0 \}$ which is the oriented angle from $(1,0)$ to $e_{1}(L)$, $\beta \in ]- \frac{\pi}{2}- \alpha, \frac{\pi}{2} - \alpha[- \{0 \}$ which is the oriented angle from $e_{1}(L)$ to $e_{2}(L)$ and $r > 0$ such that $ r^{4} < \frac{1}{\sin^{2}(\beta)}$ (see equation ($\ref{eq72}$)) which is $\lVert e_{1}(L) \rVert$. Yet we want to emphasize that all the choices of $(\alpha , \beta, r) \in \{(\alpha , \beta, r) \text{ } | \text{ }  \alpha \in ]- \frac{\pi}{2}, \frac{\pi}{2}[-\{ 0 \} \text{, } \beta \in ]- \frac{\pi}{2}- \alpha, \frac{\pi}{2} - \alpha[- \{0 \} \text{, } 0 < r \text{ and } r^{4} < \frac{1}{\sin^{2}(\beta)} \}$ does not give a lattice $L \in \mathscr{S}_{2}$ : the set of $(\alpha , \beta, r)$ that suit are a submanifold of dimension $3$.
\item For a typical lattice $L \in \mathscr{S}_{2}$ (like before), if there exists a $P \neq 0$ such that $$P(\lVert e_{1}(L) \rVert^{2}, \frac{\cos(\beta)}{\sin(\beta)}, \frac{\cos(\alpha)}{\sin(\alpha)}) = 0$$ then $L$ belongs to a countable union of submanifolds of dimension at most $2$. 
\end{itemize}
\end{proof}
We can now prove the Proposition $\ref{prop14}$.\\
\begin{proof}[Proof of Proposition $\ref{prop14}$]
We are first going to prove the first part of the proposition. \\
Let's suppose that $L$ is a typical lattice like in the Lemma $\ref{lemme13}$.
Let's prove that such a $L$ is appropriate by making an induction on the parameter $n \in \mathbb{N}$ ($r \geqslant 1$ is not fixed). \\
\\
(Basis) \\
Let's suppose that a square root $\sqrt{\prod_{i=1}^{r} A_{h_{i}}}$ belongs to $ \mathbb{Q}(X_{1},X_{2},X_{4})$. Then there exists $P,Q \in \mathbb{Q}[Y_{1},Y_{2},Y_{4}]$ such that $ P \wedge Q = 1$ and such that : 
\begin{equation}
\label{eq80}
\sqrt{\prod_{i=1}^{r} A_{h_{i}}} = (\frac{P}{Q})(X_{1},X_{2},X_{4}) \textit{.}
\end{equation}
By squaring this equation and using Lemma $\ref{lemme13}$, one gets the following equality between polynomials : 
\begin{equation}
\label{eq81}
P^{2} = Q^{2} (\prod_{i=1}^{r} A_{h_{i}})(Y_{1},Y_{2}) \textit{.}
\end{equation}
So, we get that $Q=1$. Furthermore, we have that for any $h \in \Pi-\{(1,0)\}$ : 
\begin{equation}
\label{eq82}
A_{h}(Y_{1},Y_{2})= (h_{1}Y_{1} + h_{2} Y_{2})^{2} + h_{2}^{2} > 0
\end{equation}
and this inequality is true for every $Y_{1},Y_{2} \in \mathbb{R}$. \\
So, every such $A_{h}$ is prime in the ring $\mathbb{Q}[Y_{1},Y_{2},Y_{4}]$. Furthermore, the $A_{h}$ are not associated from one another, which means that for two different $h_{1},h_{2}$, $A_{h_{1}}$ does not divide $A_{h_{2}}$. \\
The equation ($\ref{eq81}$) leads us to an absurdity.\\
\\
(Induction step) \\
Let's suppose that the result is acquired for all $k \in [0,n]$. \\
Furthermore, let's suppose that $\sqrt{\prod_{i=1}^{r} A_{h_{i}}} \in \mathbb{Q}(X_{1},X_{2},X_{4})[\sqrt{A_{k_{1}}}, \cdots, \sqrt{A_{k_{n}}},\sqrt{A_{k_{n+1}}}]$. \\
Yet, we have : $\mathbb{Q}(X_{1},X_{2},X_{4})[\sqrt{A_{k_{1}}}, \cdots, \sqrt{A_{k_{n}}},\sqrt{A_{k_{n+1}}}] = U[\sqrt{A_{k_{n+1}}}]$ where $$U=\mathbb{Q}(X_{1},X_{2},X_{4})[\sqrt{A_{k_{1}}}, \cdots, \sqrt{A_{k_{n}}}] \textit{.} $$ \\
So here there exist $x,y \in U$ such that : 
\begin{equation}
\label{eq83}
\sqrt{\prod_{i=1}^{r} A_{h_{i}}} = x + y \sqrt{A_{k_{n+1}}} \textit{.}
\end{equation}
By squaring, one gets that : 
\begin{equation}
\label{eq84}
\prod_{i=1}^{r} A_{h_{i}} = x^{2} + y^{2} A_{k_{n+1}} + 2 x y \sqrt{A_{k_{n+1}}} \textit{.}
\end{equation}
Now, if $x$ and $y$ are different from $0$, one has $\sqrt{A_{k_{n+1}}} \in U$ which is excluded by the hypothesis of the induction. \\
So $y=0$ and in that case $\sqrt{\prod_{i=1}^{r} A_{h_{i}}} \in U$, which is also excluded by the hypothesis of the induction, or $x=0$ and then $\sqrt{\prod_{i=1}^{r} A_{h_{i}}}\sqrt{A_{k_{n+1}}} \in U$, which is also excluded by the hypothesis of the induction. So, we obtain the first part of the proposition. \\
Concerning the second part, let's enumerate the $k$ that belong to $\Pi$ and that are different from $(1,0)$. It forms a sequence  $(k_{n})_{n \geqslant 1}$. Thanks to the first part of the proposition, we see that the sequence of fields $(\mathbb{Q}(X_{1},X_{2},X_{4})(\sqrt{A_{k_{1}}},\cdots,\sqrt{A_{k_{n}}}))$ is increasing for the inclusion. So, for all $n \geqslant 1$, its dimension over $\mathbb{Q}(X_{1},X_{2},X_{4})$ must be $2^{n}$. \\
Yet $\mathbb{Q}(X_{1},X_{2},X_{4})(\sqrt{A_{k_{1}}},\cdots,\sqrt{A_{k_{n}}})$ is generated by the product of square roots of $A_{k_{1}}, \cdots, A_{k_{n}}$ where $k_{i}$ taken are all different from one another. But their numbers is equal to $2^{n}$. So, they form a $\mathbb{Q}(X_{1},X_{2},X_{4})$-free family.
\end{proof}
We can now prove the Proposition $\ref{prop13}$.
\begin{proof}[Proof of Proposition $\ref{prop13}$]
Let's give us a typical $L \in \mathscr{S}_{2}$ like in the Proposition $\ref{prop14}$ and like in the Lemma $\ref{lemme13}$. Let's suppose also that : 
\begin{equation}
\sum_{i=1}^{m} p_{i} W_{k_{i}}(L) = 0
\end{equation}
with all the $p_{i} \neq 0$ and with $m \geqslant 1$.
Then the equation $(\ref{eq79})$ is verified. Yet, this equation means in particular that the family of the square roots of different $A_{k}$ is not $\mathbb{Q}(X_{1},X_{2},X_{4})$ which is excluded by Proposition $\ref{prop14}$, except in the case where $m=1$ or $m=2$ and in this last case all $k_{h}$ must be equal to $(1,0)$. So, these two cases are also excluded in the end.
\end{proof}

\subsection{Proof of Proposition $\ref{prop12}$}
Before starting the proof of Proposition $\ref{prop12}$, we only need a simple lemma : 
\begin{lemma}
\label{lemme19}
For every $m \in \mathbb{N}-\{ 0 \}$, for every family $(p_{1},\cdots,p_{m}) \in \mathbb{Z}^{m}-\{ 0 \}$, for every $k_{1},\cdots,k_{m} \in \Pi$, all distinct, for every $0 < \epsilon < 1$, there exists $a > 0$ and $K_{\epsilon}$ a measurable set of $\mathscr{S}_{2}$ such that $\tilde{\mu}_{2} (\mathscr{S}_{2}-K_{\epsilon}) \leqslant \epsilon$ and such that for all $L \in  K_{\epsilon}$,
$$ \lvert \sum_{i=1}^{m}  p_{i} W_{k_{i}}(L) \rvert \geqslant a \textit{ .} $$
\end{lemma}
\begin{proof}
It is a direct consequence of Proposition $\ref{prop13}$.
\end{proof}
Now, by using the foliation given by $\delta(\lambda)$ and previous results, we can now prove Proposition $\ref{prop12}$.
\begin{proof}[Proof of Proposition $\ref{prop12}$]
The proof in all its generality can be made as in the case where $\tilde{\mu}_{2} = \mu_{2}$. So, we will suppose for simplicity that $\tilde{\mu}_{2} = \mu_{2}$. \\
Let $l \geqslant 1$. Let $\psi \in C_{c}^{\infty}(\mathscr{S}_{2})$. Let $(p_{1},\cdots,p_{l}) \in \mathbb{Z}^{l}-\{ 0 \}$. \\
For all $\epsilon > 0$, we call $\mathcal{F}_{\epsilon}$ the tribe on $\mathscr{S}_{2}$ generated by the $S_{\epsilon}(L)$. Let $1 > \epsilon_{1} \geqslant \epsilon_{2} > 0$. \\
According to Lemmas $\ref{lemme11}$ and $\ref{lemme19}$, there exists a measurable part $K_{\epsilon_{1}}$ such that $\mu_{2}(K_{\epsilon_{1}}) \geqslant 1- \epsilon_{1}$, a real $M > 0$ and a real $a > 0$, such that
\begin{itemize}
\item for every $L \in K_{\epsilon_{1}}$, for every $\lambda \in [\frac{1}{1+\epsilon_{1}}, 1+ \epsilon_{1}]$ : 
\begin{equation}
\label{eq87}
|\psi(\delta(h)L) -\psi(L)|  \leqslant M |h| 
\end{equation}
where $h= \lambda -1$
\item for every $L \in K_{\epsilon_{1}}$, for every $\lambda \in [\frac{1}{1+\epsilon_{1}}, 1+ \epsilon_{1}]$, the equation ($\ref{eq61}$) is verified.
\item for every $L \in K_{\epsilon_{1}}$, 
\begin{equation}
\label{eq300}
\lvert \sum_{i=1}^{l}  p_{i} W_{k_{i}}(L) \rvert \geqslant a \textit{ .} 
\end{equation}
\end{itemize}
Furthermore, we are going to suppose, even if it means making $\epsilon_{2}$ goes to $0$, that for every $L \in K_{\epsilon_{1}}$, 
\begin{equation}
\label{eq506}
S_{\epsilon_{2}}(L) \subset K_{\epsilon_{1}} \textit{.}
\end{equation}
$\textbf{Claim. }$With these notations, we have, for all $L \in K_{\epsilon_{1}}$, that : 
\begin{equation}
\label{eq93}
 \mathbb{E} \left( \psi e(\sum_{i=1}^{l} p_{i} \theta_{k_{i}}) | \mathcal{F}_{\epsilon_{2}} \right) (L) = O(\epsilon_{2} ) + O(\frac{1}{a t \epsilon_{2}}) + \frac{1}{a} O(\epsilon_{2}) \textit{.}
\end{equation}
To this end,  let's set, for all $ \epsilon > 0$, $$ \delta(\epsilon) = 1+\epsilon - \frac{1}{1+\epsilon} \textit{.}$$ \\
Then one has for every $L \in K_{\epsilon_{1}}$ according to Lemma $\ref{lemme11}$ :
\begin{align}
\label{eq89}
& \mathbb{E}(\psi e(\sum_{j=1}^{l} p_{j} \theta_{k_{j}}) | \mathcal{F}_{\epsilon_{2}} )(L) \nonumber \\
& = \frac{1}{\delta(\epsilon_{2})} \int_{\frac{1}{1+\epsilon_{2}}-1}^{\epsilon_{2}} \left( \psi e(\sum_{j=1}^{l} p_{j} \theta_{k_{j}}) \right) (\delta(h) L) d h \nonumber \\
& = \psi(L) \frac{1}{\delta(\epsilon_{2})} \int_{\frac{1}{1+\epsilon_{2}}-1}^{\epsilon_{2}} \left(e(\sum_{j=1}^{l} p_{j} \theta_{k_{j}}) \right) (\delta(h) L) d h + O(\epsilon_{2}) \nonumber \\
& = \left( \psi e(\sum_{j=1}^{l} p_{j} \theta_{k_{j}}) \right)(L)  \frac{1}{\delta(\epsilon_{2})} \int_{\frac{1}{1+\epsilon_{2}}-1}^{\epsilon_{2}} e^{i t D_{1}(L) h + i t D_{2}(L,h) } d h + O(\epsilon_{2})
\end{align}
where
\begin{equation}
\label{eq90}
D_{1}(L) = \sum_{j=1}^{l} p_{j} W_{k_{j}}(L) \textit{, }
\end{equation}
\begin{equation}
\label{eq91}
D_{2}(L,h) = \sum_{j=1}^{l} p_{j} \theta_{k_{j}}(\delta(\lambda) L) -\sum_{j=1}^{l} p_{j} \theta_{k_{j}}(L)- D_{1}(L)h  \textit{ such that }
\end{equation}
\begin{equation}
\label{eq92}
D_{2}(L,h) = O(h^{2}) \textit{ and }
\end{equation}
$D_{2}(L, \cdot)$ is smooth around $0$. \\
Thus, by integrating by part and by using equation ($\ref{eq300}$), one gets that for all $L \in K_{\epsilon_{1}}$ : 
\begin{align}
\label{eq507}
& \frac{1}{\delta(\epsilon_{2})} \int_{\frac{1}{1+\epsilon_{2}}-1}^{\epsilon_{2}} e^{i t D_{1}(L) h + i t D_{2}(L,h) } d h \nonumber \\
& = \frac{1}{\delta(\epsilon_{2})} \left( \left[  \frac{e^{i t D_{1}(L) h + i t D_{2}(L,h)}}{i t D_{1}(L) } \right]_{\frac{1}{1+\epsilon_{2}}-1}^{\epsilon_{2}} + \frac{1}{D_{1}(L)} \int_{\frac{1}{1+\epsilon_{2}}-1}^{\epsilon_{2}} \left( D_{2}(L,\cdot) \right) '(h) e^{i t D_{1}(L) h + i t D_{2}(L,h)} dh \right) \nonumber \\
& = O(\frac{1}{a t \epsilon_{2}}) + \frac{1}{a} O(\epsilon_{2}) \textit{.}
\end{align}
Finally, equation ($\ref{eq89}$) and equation ($\ref{eq507}$) give the wanted claim. \\
\\
Thanks to equation ($\ref{eq93}$), the fact that $\mu_{2}(K_{\epsilon_{1}}) \geqslant 1- \epsilon_{1}$ and because of equation ($\ref{eq506}$), we have that : 
\begin{align}
\label{eq94}
& |\mathbb{E}\left(\psi e(\sum_{j=1}^{l} p_{j} \theta_{k_{j}}) \right) | \nonumber \\
& \leqslant  |\mathbb{E} \left( \psi e(\sum_{j=1}^{l} p_{j} \theta_{k_{j}}) \mathbf{1}_{K_{\epsilon_{1}}^{c}} \right)|  +  |\mathbb{E}\left( \psi e(\sum_{j=1}^{l} p_{j} \theta_{k_{j}}) \mathbf{1}_{K_{\epsilon_{1}}} \right) | \nonumber \\
& \leqslant \lVert \psi \rVert_{\infty} \epsilon_{1} + O(\epsilon_{2} ) + O(\frac{1}{a t \epsilon_{2}}) + \frac{1}{a} O(\epsilon_{2}) \textit{.}
\end{align}
 By first choosing $\epsilon_{1} > 0$ small enough (note that $a$ depends on $\epsilon_{1}$), then choosing $\epsilon_{2} > 0$ and finally choosing $t$ large enough, we obtain the wanted result. 
\end{proof}
We are now brought back to the study of the convergence, when $A \rightarrow \infty$, of a sum of the type $$\tilde{S}_{A}(\omega,L) = \sum_{k=(k_{1},k_{2}) \in \Pi_{A}(L)} \frac{Z_{k}(\omega)}{\lVert k_{1} e_{1}(L) + k_{2} e_{2}(L) \rVert^{\frac{3}{2}}}   $$
where $Z_{k}$ are real non-zero independent and identically distributed random variables from $\Omega \ni \omega$ that are symmetrical and have a compact support. We are going to see in the next section that the sums of this type converge almost surely (on the space $\Omega \times \mathscr{S}_{2}$) when $A \rightarrow \infty$. Furthermore, we are going to see that the almost sure limit $\lim_{A \rightarrow \infty} \tilde{S}_{A}(\omega,L)$ : \begin{itemize}
\item is symmetrical and its expectation is equal to $0$
\item admits moment of order $1+\kappa$ for any $0 \leqslant \kappa < \frac{1}{3}$ 
\item does not admit a moment of order $\frac{4}{3}$ when $\sigma(L) \geqslant m$ where $m > 0$ and where $L$ belongs to an event of the form $( \lVert L \rVert < \alpha )$ with $\alpha > 0$. 
\item when there exists $\alpha > 0 $ such that $\tilde{\mu}_{2}(\{ L \in \mathscr{S}_{2} \text{ } | \text{ } \lVert L \rVert_{1} < \alpha \}) = 0$ then it admits moments of all order $1 \leqslant p < \infty$.  
\end{itemize}
\section{Asymptotic study of $\tilde{S}_{A}(\omega,L)$ and proof of $\ref{thm2}$} 
The goal of this section is first to prove the following proposition : 
\begin{prop}
\label{prop15}
$\tilde{S}_{A}(\omega,L)$ converges almost surely when $A \rightarrow \infty$. Furthermore, we are going to see that the almost sure limit $\lim_{A \rightarrow \infty} \tilde{S}_{A}(\omega,L)$ : \begin{itemize}
\item is symmetrical and its expectation is equal to $0$
\item admits moment of order $1+\kappa$ for any $0 \leqslant \kappa < \frac{1}{3}$ 
\item does not admit a moment of order $\frac{4}{3}$ when $\sigma(L) \geqslant m$ where $m > 0$ and where $L$ belongs to an event of the form $( \lVert L \rVert < \alpha )$ with $\alpha > 0$. 
\item when there exists $\alpha > 0 $ such that $\tilde{\mu}_{2}(\{ L \in \mathscr{S}_{2} \text{ } | \text{ } \lVert L \rVert_{1} < \alpha \}) = 0$ then it admits moments of all order $1 \leqslant p < \infty$.  
\end{itemize} 
\end{prop}
To prove this proposition, we first need some lemmas. \\
The following lemma basically gives us that the magnitude of $\lim_{A \rightarrow \infty} \tilde{S}_{A}(\omega,L)$ is given by $\frac{X_{1,0}}{\lVert L \rVert_{1}}$.
\begin{lemma}
\label{lemme18}
Let's set, for $0 < \epsilon < 1 < A$, 
\begin{equation}
\label{eq96}
B_{k,1,A}(\omega, L, \epsilon) =  \frac{Z_{k}(\omega)}{\lVert k_{1} e_{1}(L) + k_{2} e_{2}(L) \rVert^{\frac{3}{2}}}  \mathbf{1}_{A \geqslant \lVert k_{1} e_{1}(L) + k_{2} e_{2}(L) \rVert \geqslant \epsilon } \textit{, }
\end{equation}

\begin{equation}
\label{eq99}
B_{k,2}(\omega, L, \epsilon) =  \frac{Z_{k}(\omega)}{\lVert k_{1} e_{1}(L) + k_{2} e_{2}(L) \rVert^{\frac{3}{2}}}  \mathbf{1}_{ \lVert k_{1} e_{1}(L) + k_{2} e_{2}(L) \rVert < \epsilon } \textit{, }
\end{equation}
\begin{equation} 
\label{eq104}
\tilde{S}_{\epsilon,1,A}(\omega,L)  =\sum_{k \in \Pi} H_{k,1,A}(\omega,L,\epsilon) \textit{ and}
\end{equation}
\begin{equation}
\label{eq100}
\tilde{S}_{\epsilon,2}(\omega,L)  = \sum_{k \in \Pi} H_{k,2}(\omega,L,\epsilon)
\end{equation}
so that 
$$\tilde{S}_{A}(\omega,L) = \tilde{S}_{\epsilon,1,A}(\omega,L) + \tilde{S}_{\epsilon,2}(\omega,L) \textit{.}$$
Then there exists $(M_{s})_{s \in \mathbb{N}-\{ 0,1 \}}$, a sequence of positive real numbers (that does not depend on $A$), such that for all $s \geqslant 2 $: 
\begin{equation}
\label{eq98}
\mathbb{E}(\sum_{k \in \Pi} |H_{k,1}|^{s}(\omega, L,A, \epsilon)) \leqslant M_{s} \textit{.}
\end{equation}
Furthermore, one has : 
\begin{equation}
\label{eq103}
\tilde{S}_{\epsilon,2}(\omega,L) =  \frac{Z_{1,0}}{\lVert L \rVert_{1}^{\frac{3}{2}}} \mathbf{1}_{\lVert L \rVert_{1} < \epsilon}
\end{equation}
when $\epsilon > 0$ is chosen small enough.
\end{lemma}
To prove this lemma, we need to recall a Theorem due to Minkowski. 
\begin{theorem}[\cite{minkowski2016geometrie}]
\label{thm3}
There exists $K > 0$  such that for all $L \in \mathscr{S}_{2}$, 
$$\lVert L \rVert_{1} \lVert L \rVert_{2}  \geqslant K_{2} \textit{.}$$
\end{theorem}
We can now prove Lemma $\ref{lemme18}$.
\begin{proof}[Proof of Lemma $\ref{lemme18}$]
The first fact is a direct consequence of Lemma $\ref{lemme5}$ and of the fact that the $Z_{k}$ are of compact support. \\
Indeed, first let's say that it is enough to deal with the case $\tilde{\mu}_{2}= \mu_{2}$. With the notations of Lemma $\ref{lemme5}$, let's set 
\begin{equation}
\label{eq203}
f(x) = \frac{1}{\lVert x \rVert^{3}} \mathbf{1}_{\epsilon \leqslant \lVert x \rVert \leqslant A }  
\end{equation}
and
\begin{equation}
\label{eq204}
\mathcal{S}(f)(L) = \sum_{l \text{ prime  } \in L} f(l) \textit{.}
\end{equation}
From the facts that the $X_{k}$ are independent, identically distributed, symmetrical with a compact support, one gets that : 
\begin{equation}
\label{eq205}
\mathbb{E}(\sum_{k \in \Pi} (B_{k,1,A})^{2}(\omega, L, \epsilon)) = \mathbb{E}(Z_{1,0}^{2}) \mathbb{E}(F) \textit{.}
\end{equation}
The Lemma $\ref{lemme5}$ gives then that : 
\begin{equation}
\label{eq206}
\mathbb{E}(\sum_{k \in \Pi} (B_{k,1,A})^{2}(\omega, L, \epsilon)) = \mathbb{E}(Z_{1,0}^{2}) \int_{\mathbb{R}^{2}} f(x) dx \textit{.}
\end{equation}
By passing into polar coordinates, one gets finally that : 
\begin{equation}
\label{eq207}
\mathbb{E}(\sum_{k \in \Pi} (B_{k,1,A})^{2}(\omega, L, \epsilon)) = 2 \pi \mathbb{E}(Z_{1,0}^{2}) ( \frac{1}{\epsilon} - \frac{1}{A} ) \textit{.}
\end{equation}
We do the same for the other $k \geqslant 3$. \\
Concerning the second fact, if there exists $(k_{1},k_{2}) \in \Pi $ such that $\lVert k_{1} e_{1}(L) + k_{2} e_{2}(L) \rVert < \epsilon$ then $\lVert L \rVert_{1} < \epsilon$ and the reverse is also true. \\
In this case, according to Theorem $\ref{thm3}$, one has 
\begin{equation}
\label{eq102}
\lVert k_{1}e_{1}(L) + k_{2} e_{2}(L) \rVert \geqslant \lVert L \rVert_{2} > \frac{K_{2}}{\epsilon}
\end{equation}
where $(k_{1},k_{2}) \in \Pi-\{(1,0)\}$. \\
Thus, by choosing $\epsilon > 0$ so that $\frac{K_{2}}{\epsilon}> \epsilon$, we get that : 
\begin{equation}
\label{eq103}
\tilde{S}_{ \epsilon,2}(\omega,L) =  \frac{Z_{1,0}}{\lVert L \rVert_{1}^{\frac{3}{2}}}
\end{equation} 
when $\lVert L \rVert_{1} < \epsilon$. \\
Otherwise, $\tilde{S}_{ \epsilon,2}(\omega,L)= 0 $ and $\frac{Z_{1,0}}{\lVert L \rVert^{\frac{3}{2}}} \mathbf{1}_{\lVert L \rVert_{1} < \epsilon} = 0$. And so the second result is also true.
\end{proof}
The following lemma precise the Lemma $\ref{lemme8}$ and it gives us a better understanding of the law of $\frac{X_{1,0}}{\lVert L \rVert_{1}^{\frac{3}{2}}}$. The proof of it is basically the same as the proof of Lemma $\ref{lemme8}$. 
\begin{lemma} 
\label{lemme14}
When $\sigma(L) \geqslant m$ where $m > 0$ and where $L$ belongs to an event of the form $( \lVert L \rVert < \alpha )$ with $\alpha > 0$, there exists $C > 0$  such that
$$\mathbb{P}( \lVert L \rVert_{1}^{-\frac{3}{2}} > \beta) = \frac{C}{\beta^{\frac{4}{3}}} $$
and this equality is true for every $\beta $ large enough.
\end{lemma}
\begin{proof}
The proof is a direct application of Lemma $\ref{lemme5}$.
\end{proof}
We recall a classic tool from probability theory : 
\begin{lemma}
\label{lemme15}
For $X$ a real random variable, one has for every $k \geqslant 1$
$$ \mathbb{E}(|X|^{k}) = \int_{0}^{\infty} t^{k-1} P(|X| > t) dt \textit{.}$$
\end{lemma}
Thus, by using the fact that $Z_{1,0}$ is bounded and is different from $0$ and by using Lemma $\ref{lemme14}$ and Lemma $\ref{lemme15}$, we get the following fact about the behaviour $\frac{Z_{1,0}}{\lVert L \rVert_{1}^{\frac{3}{2}}}$ : 
\begin{lemma}
\label{lemme17}
For every $0 \leqslant \kappa < \frac{1}{3}$, $\frac{Z_{1,0}(\omega)}{\lVert L \rVert_{1}^{\frac{3}{2}}}$ admits a moment of order $1+\kappa$. \\
Furthermore, $\frac{Z_{1,0}(\omega)}{\lVert L \rVert_{1}^{\frac{3}{2}}}$ does not admit a moment of order $\frac{4}{3}$ when $\sigma(L) \geqslant m$ where $m > 0$ and where $L$ belongs to an event of the form $( \lVert L \rVert < \alpha )$ with $\alpha > 0$. \\
\end{lemma}
We finally remind the reader of the following obvious lemma : 
\begin{lemma}
\label{lemme16}
If $X$ is a real random variable integrable then if it is symmetrical, its expectation is equal to zero.
\end{lemma}
We can now prove the Proposition $\ref{prop15}$.
\begin{proof}[Proof of Proposition $\ref{prop15}$]
Because of the Lemma $\ref{lemme16}$, we only need first to prove that $\tilde{S}_{A}(\omega,L)$ converges almost surely when $A \rightarrow \infty$ and that the limit random variable is symmetrical and admits moment of order $1+\alpha$ and that for all $0 \leqslant \alpha < \frac{1}{3}$ : the fact that the expectation of the limit is null will be a consequence of Lemma $\ref{lemme16}$. \\
\\
$\bullet$ Let's prove that $\tilde{S}_{A}(\omega,L)$ converges almost surely when $A \rightarrow \infty$. It is enough to prove that it is the case on all the events $\{ \lVert L \rVert_{1} \geqslant \epsilon \}$ where $\epsilon > 0$. So, let's give us $\epsilon > 0$. \\
On such an event, one has : 
\begin{equation}
\label{eq95}
\tilde{S}_{A}(\omega,L) = \tilde{S}_{\epsilon,1,A}(\omega,L)  \textit{.}
\end{equation}
Yet, one has, because the $X_{k}$ are symmetrical and integrable,
\begin{equation}
\label{eq97}
\mathbb{E}(B_{k,1,A}(\omega, L, \epsilon)) = 0 \textit{.}
\end{equation}
By using ($\ref{eq98}$) (see Lemma $\ref{lemme18}$) and the equation $(\ref{eq97})$ and because the $Z_{k}$ are independent between them and their expectations are equal to zero, one gets that :
$$\tilde{S}_{\epsilon,1,A}(\omega,L) \textit{ converges almost surely when } A \rightarrow \infty \textit{.}$$
Because of equation $(\ref{eq95})$, one gets that $\tilde{S}_{A}(\omega,L)$ converges almost surely when $A \rightarrow \infty$ on all the events $\{ \lVert L \rVert_{1} \geqslant \epsilon \}$ where $\epsilon > 0$. 
\\
\\
$\bullet$ The fact that $\lim_{A \rightarrow \infty} \tilde{S}_{A}(\omega,L)$ is symmetrical follows directly from the fact that the $Z_{k}$ are symmetrical. \\
\\
$\bullet$ It only remains to prove that $\lim_{A \rightarrow \infty} \tilde{S}_{A}(\omega,L)$ admits a moment of order $1+\kappa$ and that for all $ 0 \leqslant \kappa < \frac{1}{3}$ ; that $\lim_{A \rightarrow \infty} \tilde{S}_{A}(\omega,L)$ does not admit a moment of order $\frac{4}{3}$ when $\sigma(L) \geqslant m$ where $m > 0$ and where $L$ belongs to an event of the form $( \lVert L \rVert < \alpha )$ with $\alpha > 0$.  ; and that if $\tilde{\mu}_{2}$ is such that there exists $\alpha > 0$ such that $\tilde{\mu}_{2}(\{ L \in \mathscr{S}_{2} \text{ } | \text{ } \lVert L \rVert_{1} < \alpha \} ) = 0$ then the limit  $\lim_{A \rightarrow \infty} \tilde{S}_{A}(\omega,L)$ admits moments for all orders $1 \leqslant p < \infty$. These facts are a direct consequence of Lemma $\ref{lemme18}$, of Lemma $\ref{lemme17}$ and of Fatou's lemma.
\end{proof}

\bibliographystyle{plain}
\bibliography{bibliographie}
\end{document}